\documentclass[review]{elsarticle}

\usepackage{hyperref}

\usepackage[tbtags]{amsmath}
\usepackage{amsfonts,amssymb,mathrsfs,amscd,comment}

\newtheorem{Lm}{Lemma}[section]

\newtheorem{Th}{Theorem}[section]

\newtheorem{Remark}{Remark}[section]

\newtheorem{Def}{Definition}[section]

\begin{document}

\begin{frontmatter}

\title{Besov-type spaces of variable smoothness on rough domains}
\tnotetext[mytitlenote]{
This work is supported by the Russian Science Foundation under grant 14-50-00005. }

\author{A.\,I.~Tyulenev}
\address{Steklov Mathematical Institute {\rm (}Russian Academy of Sciences{\rm )}}

\ead{tyulenev-math@yandex.ru, tyulenev@mi.ras.ru}

\numberwithin{equation}{section}

\begin{abstract}
The paper puts forward new Besov spaces of variable smoothness
$B^{\varphi_{0}}_{p,q}(G,\{t_{k}\})$ and $\widetilde{B}^{l}_{p,q,r}(\Omega,\{t_{k}\})$ on rough domains.
A~domain~$G$ is either a~bounded Lipschitz domain in~$\mathbb{R}^{n}$ or the epigraph of a~Lipschitz function,
a~domain~$\Omega$ is an $(\varepsilon,\delta)$-domain. These spaces are shown to be the traces of the spaces
$B^{\varphi_{0}}_{p,q}(\mathbb{R}^{n},\{t_{k}\})$ and  $\widetilde{B}^{l}_{p,q,r}(\mathbb{R}^{n},\{t_{k}\})$ on domains $G$ and~$\Omega$, respectively.
The extension operator $\operatorname{Ext}_{1}:B^{\varphi_{0}}_{p,q}(G,\{t_{k}\}) \to B^{\varphi_{0}}_{p,q}(\mathbb{R}^{n},\{t_{k}\})$ is linear, the operator $\operatorname{Ext}_{2}:\widetilde{B}^{l}_{p,q,r}(\Omega,\{t_{k}\}) \to \widetilde{B}^{l}_{p,q,r}(\mathbb{R}^{n},\{t_{k}\})$ is nonlinear.
As a~corollary, an exact description of the traces of 2-microlocal Besov-type spaces and weighted Besov-type spaces on rough domains is obtained.
\end{abstract}


\end{frontmatter}

\section{Introduction}
Besov spaces of variable smoothness (and various extensions thereof) have been intensively studied in the last 30 years.
Among the fundamental works we mention \cite{Almeida}, \cite{Kempka}, \cite{Kempka2},  \cite{Kempka3}, \cite{Kempka4}, \cite{Kempka5}, \cite{Drihem}, \cite{Drihem2}, \cite{Drihem3}, \cite{Drihem4}, \cite{Drihem5}, \cite{Be1}, \cite{Be2}, \cite{Ty1}, \cite{Ty2}, \cite{Dachun}, \cite{Dachun2}, and~\cite{Sawano}.
We also note the recent paper~\cite{Ty3}, which introduces a~new Besov-type space, which has not only variables smoothness, but also in a~sense a~`variable structure'.
This space does not agree with any of the previously introduced Besov spaces of variable smoothness.
Despite of the vast literature on spaces of variable smoothness, the question of the description of traces of such spaces on rough domains remained open.
At the same time (as far as the author is aware) no question was given in the problem of exact description of traces of weighted Besov and Lizorkin--Triebel spaces
(in the entire range of parameters of smoothness and integration) with weights locally satisfying the Muckenhoupt condition on rough domains.
The author knows here only the papers \cite{Chua}, \cite{Chua2}, in which the similar problem was solved for weighted Sobolev spaces.

It is interesting to note that even for the classical Besov and Lizorkin--Triebel spaces the problem of exact description of the trace on nonsmooth domains
(in the entire range of parameters $p,q$ and $s$) was solved only recently in \cite{DeVore}, \cite{Se}, \cite{Ry2}.

We recall (see \cite{Tri}, Ch.~2 for details) that there exist several nonequivalent (in general) approaches to the construction of the theory of classical Besov spaces.
The first approach, which depends on the Fourier analysis, is convenient in working with distributions from the space $S'(\mathbb{R}^{n})$.
The second approach is based on the  approximation theory and is suitable for dealing with functions that are locally integrable in some power.
The majority of the available studies on spaces of variable smoothness was based on the Fourier analytic approach,
the only exception are the papers \cite{Kempka}, \cite{Be1}, \cite{Be2}, in which the classical approach was employed.

Recently, the author \cite{Ty1} proposed a~new nonlinear approximation approach to the study of Besov spaces of variable smoothness.
In~\cite{Ty2} a~new approach involving convolutions with smooth kernels (instead of the Fourier transform) was proposed to determine the norm in the
corresponding space. The Besov spaces of~\cite{Ty2} extend (in the case of constant exponents $p,q$)
both the well-known 2-microlocal Besov spaces, which were first studied by H.~Kempka in~\cite{Kempka4}, and the weighted Besov spaces from~\cite{Ry}.
Besides,  a~more subtle method was introduced in \cite{Ty1}, \cite{Ty2} for the analysis of a~weight sequence $\{t_{k}\}$ (which specifies the variable smoothness).
This methods was found to be instrumental in refining many available theorems on Besov spaces of variable smoothness by relaxing unnatural pointwise constraints on a~variable smoothness $\{t_{k}\}$.

\smallskip

The purpose of the present paper is to provide an exact description of the traces of Besov spaces of variable smoothness on domains with  nonsmooth boundary.
For the spaces $B^{\varphi_{0}}_{p,q}(\mathbb{R}^{n},\{t_{k}\})$, which were first introduced in~\cite{Ty2}, we put forward a~description of the traces
on bounded Lipschitz domains and on special Lipschitz domains (epigraphs of Lipschitz functions).
An exact description of the traces on $(\varepsilon,\delta)$-domains will be given for the
spaces  $\widetilde{B}^{l}_{p,q,r}(\mathbb{R}^{n},\{t_{k}\})$, which were introduced by the author in~\cite{Ty1}.

Since the weighted Besov spaces $B^{s}_{p,q}(\mathbb{R}^{n},\gamma)$, $\widetilde{B}^{s}_{p,q,r}(\mathbb{R}^{n},\gamma)$
 with a~weight~$\gamma$ locally satisfying the corresponding Muckenhoupt  condition, as well as the 2-microlocal spaces $B^{\{s_{k}\}}_{p,q}(\mathbb{R}^{n})$
 (wich were first introduced by O.~Besov \cite{Be1} for $p,q \in (1,\infty)$ and then by H.~Kempka \cite{Kempka4} for $p,q \in (0,\infty]$)
 are particular cases of the above spaces of variable smoothness (the proof of this fact can be found in \cite{Ty1}, \cite{Ty2} and in Remarks \ref{Remark2.6}, \ref{Rem3.4} below),
 we obtain as a~straightforward corollary an exact description of the traces of weighted Besov spaces and 2-microlocal Besov-type spaces on domains with nonsmooth boundary.

\section{Notation, definitions and statements of the problems}

The symbols $\mathbb{N}$, $\mathbb{N}_{0}$, $\mathbb{R}^{n}$, $\mathbb{Z}^{n}$ have their standard meaning.
By $Q$~we shall denote a~\textit{closed} cube in the space $\mathbb{R}^{n}$ with sides parallel
to the coordinate axes; $Q_{k,m}$ will denote a~\textit{closed dyadic} cube of rank $k \in \mathbb{N}_{0}$ in $\mathbb{R}^{n}$. More precisely, we put $Q_{k,m}:=\prod\limits_{i=1}^{n}\left[\frac{m_{i}}{2^{k}},\frac{m_{i}+1}{2^{k}}\right]$ for $k \in \mathbb{N}_{0}$, $m \in \mathbb{Z}^{n}$. We also set $\widetilde{Q}_{k,m}:=\prod\limits_{i=1}^{n}\left[\right.\frac{m_{i}}{2^{k}},\frac{m_{i}+1}{2^{k}}\left.\right)$ for $k \in \mathbb{N}_{0}$, $m \in \mathbb{Z}^{n}$.
In what follows, $I$ will always denote the cube $[-1,1]^{n}$.
Domains in~$\mathbb{R}^{n}$ will be denoted by the letters $G$ or~$\Omega$.

Next, the symbols $k,j,l,i$ will denote integer variables, $m$ will denote vectors from $\mathbb{Z}^{n}$.

Given a $p \in (0,\infty]$ and a~measurable set $E \subset \mathbb{R}^{n}$ of positive measure, we let $L_{p}(E)$ denote set of all classes of equivalent measurable functions with finite quasi-norm $\|f|L_{p}(E)\|:=\Bigl(\int\limits_{E}|f(x)|^{p}\,dx\Bigr)^{\frac{1}{p}}$ (the `ess\,sup' modifications in the case $p=\infty$ are clear).
For $q \in (0,\infty]$, by $l_{q}$ we shall denote the linear space of all real sequences $\{a_{k}\}$ with finite quasi-norm
$\|a_{k}|l_{q}\|:=\Bigl(\sum\limits_{k=1}^{\infty}|a_{k}|^{q}\Bigr)^{\frac{1}{q}}$ (the `sup'-modifications in the case $q=\infty$ are standard).

We will use the following convention throughout this paper.
Integrations are carried out over the whole $\mathbb{R}^{n}$, unless other limits are indicated.
Similarly, when considering some function space on $\mathbb{R}^{n}$, for brevity we shall drop the symbol $\mathbb{R}^{n}$ in the notation
of this space. For example,  instead of $C^{\infty}(\mathbb{R}^{n})$, $L_{p}(\mathbb{R}^{n})$ and so on, we shall write $C^{\infty}$, $L_{p}$, etc., respectively.

Given a function $g:\mathbb{R}^{n} \mapsto \mathbb{R}$, $j \in \mathbb{N}$, we set $g_{j}:=2^{jn}g(\cdot 2^{j})$.

Also, for a function $\psi \in C^{\infty}$ by $L_{\psi}$ we shall denote the supremum of the numbers~$L$ for which
\begin{equation}
\label{eq2.1}
\int x^{\beta}\psi(x)\,dx=0 \mbox{ for } |\beta| \le L.
\end{equation}

A~function $\psi$ has zero moments up to the order $L$ if condition \eqref{eq2.1} holds. If a~function~$\psi$ has no zero moments, then we write $L_{\psi}=-1$.

By a weight we shall understand an arbitrary measurable function which is positive almost everywhere. The basic definitions and properties of
the weighted class $A^{\mathrm{loc}}_{p}(\mathbb{R}^{n})$
with $p \in (1,\infty]$ may be found in~\cite{Ry}.

The symbols $c$ or $C$ will be used to denote (in general different) `insignificant’ constants in various inequalities.
We shall not label different constants with different indexes. When required we shall indicate the parameters on which some or
other constant depends.

The symbols $S$, $D(\Omega)$, $S'$, $D'(\Omega)$ ($\Omega$ is an open set) will have the standard meaning and denote the linear spaces of test functions and
the dual spaces of distributions (see \cite{Rudin}, Chapters 6 and~7 for details).

For further purposes we introduce the following special class of weight sequences. By a~weight sequence (weight sequences will be denoted by $\{s_{k}\}$, $\{t_{k}\}$)
we shall understand a~function sequence in which any element is a~measurable function on $\mathbb{R}^{n}$ which is positive almost everywhere.

\begin{Def}\rm (\cite{Kempka4})
\label{Def2.1}
A weight sequence $\{s_{k}\}=\{s_{k}(\cdot)\}_{k=0}^{\infty}$ will be said to lie in $\mathcal W^{\alpha_{3}}_{\alpha_{1},\alpha_{2}}$ if,
for $\alpha_{3}\geq 0$, $\alpha_{1},\alpha_{2} \in \mathbb{R}$,
\begin{equation}
\label{eq2.2}
\begin{aligned}
1)&\qquad \frac{1}{C_{1}} 2^{\alpha_{1}(k-l)} \le \frac{s_{k}(x)}{s_{l}(x)} \le C_{1} 2^{\alpha_{2}(k-l)},\ \  l \le k \in \mathbb{N}_{0},\ \  x \in \mathbb{R}^{n};\\
2)&\qquad  s_{k}(x) \le C_{2} s_{k}(y)(1+2^{k}|x-y|)^{\alpha_{3}},\ \ k \in \mathbb{N}_{0},\ \  x,y \in \mathbb{R}^{n}.
\end{aligned}
\end{equation}
\end{Def}

The constants $C_{1},C_{2}$ in \eqref{eq2.2} are independent of $k,l$ and $x,y$.

In the majority of presently available studies, Besov spaces of variable smoothness were defined in terms of the Fourier transform.
The starting points for the Fourier-analytic approach to Besov spaces of variable smoothness were the papers \cite{Be1} and~\cite{Kempka4}.

In what follows we shall need the standard resolution of unity. More precisely, let $B$ be the unit ball in $\mathbb{R}^{n}$, $\Psi_{0} \in S$ and let $\Psi_{0}(x)=1$ for $x \in B$, $\hbox{suppp}\Psi_{0} \subset 2B$. Given $j \in \mathbb{N}$, we set $\Psi_{j}(x)=\Psi_{0}(2^{-j}x) - \Psi_{0}(2^{-j+1}x)$ for $x \in \mathbb{R}^{n}$.

\begin{Def}
\label{Def2.2} \rm
Let $p,q \in (0,\infty]$, $\alpha_{1},\alpha_{2} \in \mathbb{R}$, $\alpha_{3} \geq 0$, $\{s_{k}\} \in \mathcal W^{\alpha_{3}}_{\alpha_{1},\alpha_{2}}$.
By $B^{\{s_{k}\}}_{p,q}$ we shall denote the space of all distributions $f \in S'$ with finite quasi-norm
\begin{equation}
\label{eq2.3}
\|f|B^{\{s_{k}\}}_{p,q}\|:=\|s_{k}F^{-1}(\Psi_{k}F[f])|l_{q}(L_{p})\|.
\end{equation}
\end{Def}

In \eqref{eq2.3} the symbols $F$ and $F^{-1}$ denote the direct and inverse Fourier transforms, respectively.

However, Definition \ref{Def2.2} is not satisfactory for the following reasons. First, the weight sequence $\{s_{k}\}$ contains functions that
grow slowly at infinity. Besides, any function $s_{k}$ has no singular points.

To be able to work with more involved weight sequences, a~definition of a new weighted class (see Definition \ref{Def2.3} below)
$\mathcal{X}^{\alpha_{3}}_{\alpha,\sigma,p}$ was proposed in \cite{Ty1}, \cite{Ty2}; this class can be looked upon as a~multi-weighed generalization
of the local Muckenhoupt class $A^{\mathrm{loc}}_{p}(\mathbb{R}^{n})$ (which was introduced in \cite{Ry}). The study of Besov spaces with variable smoothness from the class $\mathcal{X}^{\alpha_{3}}_{\alpha,\sigma,p}$
is instrumental for developing a~unified approach both to 2-microlocal Besov spaces and to weighted Besov-type spaces.

If given $p \in (0,\infty]$ a~sequence of weights $\{t_{k}\}$  is such that $t_{k} \in L^{\mathrm{loc}}_{p}$ for  $k \in \mathbb{N}_{0}$,
then the weight sequence $\{t_{k}\}$ will be called a~$p$-admissible weight sequence.

Let $\{t_{k}\}$ be a $p$-admissible weight sequence. By $\{t_{k,m}\}$ we shall denote the multiple sequence defined by

$$
t_{k,m} := \|t_{k}|L_{p}(Q_{k,m})\|  \quad \hbox{ for } k \in \mathbb{N}_{0}, m \in \mathbb{Z}^{n};
$$

by $\{\overline{t}_{k}\}$ we shall denote a weight sequence of the form

$$
\overline{t}_{k}(x):=2^{\frac{kn}{p}}\sum\limits_{m \in \mathbb{Z}^{n}}t_{k,m}\chi_{\widetilde{Q}_{k,m}}(x)
  \quad \hbox{ for } k \in \mathbb{N}_{0}, x \in \mathbb{R}^{n}.
$$

\begin{Def}
\label{Def2.3} \rm
Let $p,\sigma_{1},\sigma_{2} \in (0,\infty]$, $\alpha_{1}, \alpha_{2} \in \mathbb{R}$, $\alpha_{3} \geq 0$, and let
$\sigma:=(\sigma_{1},\sigma_{2})$, $\alpha:=( \alpha_{1},\alpha_{2})$. We let
$\mathcal X^{\alpha_{3}}_{\alpha,\sigma,p}$ denote the set of all $p$-admissible weight sequences $\{t_{k}\}:=\{t_{k}(\cdot)\}_{k=0}^{\infty}$ such that, for some $C_{1},C_{2}>0$,

{\rm 1)} \begin{equation}
\label{eq2.4}
\Bigl(2^{kn}\int\limits_{Q_{k,m}}t^{p}_{k}(x)\Bigr)^{\frac{1}{p}}\Bigl(2^{kn}\int\limits_{Q_{k,m}}t^{-\sigma_{1}}_{j}(x)\Bigr)^{\frac{1}{\sigma_{1}}} \le C_{1}2^{\alpha_{1}(k-j)}, \qquad
 0 \le k \le j, m \in \mathbb{Z}^{n},
\end{equation}

{\rm 2)} \begin{equation}
\label{eq2.5}
\Bigl(2^{kn}\int\limits_{Q_{k,m}}t^{p}_{k}(x)\Bigr)^{-\frac{1}{p}}\Bigl(2^{kn}\int\limits_{Q_{k,m}}t^{\sigma_{2}}_{j}(x)\Bigr)^{\frac{1}{\sigma_{2}}} \le C_{2}2^{\alpha_{2}(j-k)}, \qquad
 0 \le k \le j, m \in \mathbb{Z}^{n},
\end{equation}
(the modifications of \eqref{eq2.4} or \eqref{eq2.5} in the case $\sigma_{1}=\infty$ or $\sigma_{2}=\infty$ are standard).

{\rm 3)} for all $k \in \mathbb{N}_{0}$
 \begin{equation}
\label{eq2.6}
0<t_{k,m} \le 2^{\alpha_{3}} t_{k,\widetilde{m}}, \qquad m,\widetilde{m} \in \mathbb{Z}^{n}, \ \ |m_{i}-\widetilde{m}_{i}| \le 1, \ \ i=1,\dots,n,
\end{equation}
\end{Def}

\begin{Remark} \rm
\label{Remark2.1}
From the H\"older inequality it follows that
$\mathcal{X}^{\alpha_{3}}_{\alpha,\widetilde{\sigma},p} \subset \mathcal{X}^{\alpha_{3}}_{\alpha,\sigma,p}$ if $\widetilde{\sigma}_{i} \geq \sigma_{i}$, $i=1,2$.

We shall note the estimate
\begin{equation}
\label{eq2.7}
t_{k,m} \Bigl(\sum\limits_{\substack{\widetilde{m} \in \mathbb{Z}^{n} \\ Q_{j,\widetilde{m}} \subset Q_{k,m}}}\frac{1}{t^{\sigma_{1}}_{j,\widetilde{m}}}\Bigr)^{\frac{1}{\sigma_{1}}} \le C_{1}2^{(\alpha_{1}-\frac{n}{p}-\frac{n}{\sigma_{1}})(k-j)}, \quad 0 \le k \le j, m \in \mathbb{Z}^{n},
\end{equation}
which is a direct corollary to \eqref{eq2.4}, provided that  $\{t_{k}\}=\{\overline{t}_{k}\} \in \mathcal{X}^{\alpha_{3}}_{\alpha,\sigma,p}$
are  $p$-admissible sequences.

The following inequality is an easy consequence of \eqref{eq2.6} and \eqref{eq2.7}:
\begin{equation}
\label{eq2.8}
\frac{1}{C}t_{j,\widetilde{m}} \le t_{k,m} \le Ct_{j,\widetilde{m}}, \quad j \geq k;
\end{equation}
here  $\widetilde{m} \in \mathbb{Z}^{n}$ is an arbitrary index, for which $Q_{j,\widetilde{m}} \subset Q_{k,m}$.
The constant $C > 0$ in~\eqref{eq2.8} depends only on $n$, $p$, $\sigma$, $\alpha$, $\alpha_{3}$ and (what is important) on~$j-k$.
\end{Remark}

\begin{Remark} \rm
\label{Remark2.2}
Let $\{t_{k}\} \in \mathcal{X}^{\alpha_{3}}_{\alpha,\sigma,p}$. The estimate
\begin{equation}
\label{eq2.9}
t_{k,m} \le C \hbox{exp}(\alpha_{3} \hbox{ln} 2 |m-\widetilde{m}|)t_{k,\widetilde{m}}, \quad m,\widetilde{m} \in \mathbb{Z}^{n}
\end{equation}
follows from  condition \eqref{eq2.6} by induction.
\end{Remark}

\begin{Remark} \rm
\label{Remark2.2'}
It is easily seen that $\mathcal{W}^{\alpha_{3}}_{\alpha_{1},\alpha_{2}} \subset \mathcal{X}^{\alpha_{3}}_{\alpha,\sigma,p}$ for any choice of
 parameters $\alpha_{3},\alpha,\sigma,p$. The converse inclusion does not hold, as shows the example of a~weight sequence $\{t_{k}\}$=$2^{ks}\gamma$, where
 $s \in \mathbb{R}$, $\gamma^{p} \in A^{loc}_{\nu}$ with some $\nu \in (1,\infty)$. Besides, for the sequences $\{t_{k}\}=\{\overline{t}_{k}\}$ we get the equivalence of \eqref{eq2.5}, \eqref{eq2.6} with $\sigma_{1},\sigma_{2} = \infty$ and the first condition in \eqref{eq2.2}.
 $\mathcal{W}^{\alpha_{3}}_{\alpha_{1},\alpha_{2}} \subset \mathcal{X}^{\alpha_{3}}_{\alpha,\sigma,p}$ with $\sigma_{1}=\sigma_{2}=\infty$.
 So, the class $\mathcal{X}^{\alpha_{3}}_{\alpha,\sigma,p}$ is a~natural extension of the class $\mathcal{W}^{\alpha_{3}}_{\alpha_{1},\alpha_{2}}$.
\end{Remark}

Following \cite{Sc} by $S_{e}$ we denote the set of all $f \in C^{\infty}$ such that
\begin{equation}
\label{eq2.10}
p_{N}(f):=\sup\limits_{x \in \mathbb{R}^{n}}\hbox{exp}(N|x|)\sum\limits_{|\alpha| \le N}|D^{\alpha}f(x)| < \infty \quad \hbox{for all } N \in \mathbb{N}_{0}.
\end{equation}

We equip $S_{e}$ with the locally convex topology defined by the system of the semi-norms $p_{N}$.

By  $S'_{e}$ we denote the collections of all continuous linear forms on $S_{e}$. We equip $S'_{e}$ with the strong topology (see \cite{Sc} for details).

\begin{Def} \rm
\label{Def2.4}
Let $p,q,\sigma_{1},\sigma_{2} \in (0,\infty]$, $\alpha_{1},\alpha_{2} \in \mathbb{R}$, $\alpha_{3} \geq 0$, and let $\{t_{k}\} \in \mathcal X^{\alpha_{3}}_{\alpha,\sigma,p}$
be a~weight sequence. Consider $\varphi_{0} \in D$ such that $\int \varphi_{0}(x)\,dx =1$. Next, we set
$\varphi(x):=\varphi_{0}(x)-2^{-n}\varphi_{0}(\frac{x}{2})$, where $x \in \mathbb{R}^{n}$. By $B^{\varphi_{0}}_{p,q}(\{t_{k}\})$ we shall denote the set of all
distributions  $f \in S'_{e}$ with finite quasi-norm
\begin{equation}
\label{eq2.11}
\|f|B^{\varphi_{0}}_{p,q}(\{t_{k}\})\|:= \left(\sum\limits_{k=0}^{\infty} \|t_{k} (\varphi_{k} \ast f)|L_{p}\|^{q}\right)^{\frac{1}{q}}
\end{equation}
(the modifications of \eqref{eq2.11} in the case $q=\infty$ are straightforward).
\end{Def}

\begin{Remark} \rm
\label{Remark2.3}
The question naturally arises of the well-definedness of the space $B^{\varphi_{0}}_{p,q}(\{t_{k}\})$ from Definition \ref{Def2.4}; that is, whether the space $B^{\varphi_{0}}_{p,q}(\{t_{k}\})$ is independent of the choice of the function
$\varphi_{0}$ and whether the corresponding quasi-norms are equivalent. From Theorem~3.2 of~\cite{Ty2} and Remark~\ref{Remark2.1}
it follows that the space $B^{\varphi_{0}}_{p,q}(\{t_{k}\})$  is well-defined in the above sense for $1+L_{\varphi} > \alpha_{2}$ and $\sigma_{2} \geq p$.
\end{Remark}

Let  $Q$ be a cube of side length $l(Q)=2^{-k}$ with some $k \in \mathbb{N}_{0}$. Given a~fixed $r \in (0,\infty]$, we let $P_{Q}[f]$ denote a~polynomial of near-best (with some constant $\lambda \geq 1$) approximation
to a~function~$f$ in the quasi-norm $L_{r}(Q)$  (see \cite{DeVore} for details). Next, let $E_{l}(f,Q)_{r}$ be the best approximation to function $f \in L^{\mathrm{loc}}_{r}$ on a~cube~$Q$ in the
quasi-norm $L_{r}(Q)$. We set $\mathcal{E}_{l}(f,Q)_{r}:=2^{\frac{kn}{r}} E_{l}(f,Q)_{r}$.

Let $U \subset \mathbb{R}^{n}$ be an open set, $l \in \mathbb{N}$. Given $x,h \in \mathbb{R}^{n}$, we set
$$
\Delta^{l}(h,U)f(x):=
\begin{cases}
\displaystyle \sum\limits_{i=0}^{l}C^{i}_{l}f(x+ih) &\mbox{ if}\ \  [x,x+hl] \subset U, \\
0 & \mbox{otherwise}
\end{cases}
$$

Let $Q \subset U$ be a cube of side length $l(Q)=2^{-k}$ for some $k \in \mathbb{N}_{0}$. Given $r \in (0,\infty]$, we set (the modifications in the case $r=\infty$
are clear)
$$
\delta^{l}_{r}(Q,U) f:=\Bigl(2^{2kn}\int\limits_{Q}\int\limits_{2^{-k}I}|\Delta^{l}(h,U)f(x)|^{r}\,dh\,dx\Bigr)^{\frac{1}{r}}.
$$

In the case $U=\mathbb{R}^{n}$ we shall write $\delta^{l}_{r}(Q)f$ for $\delta^{l}_{r}(Q,\mathbb{R}^{n}) f$.

Let $f \in L^{\mathrm{loc}}_{r}$. The following fundamental estimates will be of great value in the future.

\begin{equation}
\label{eq2.12}
\mathcal{E}_{l}(f, cQ_{0,m})_{r} \le C \sum\limits_{\substack{\widetilde{m} \in \mathbb{Z}^{n} \\Q_{k,\widetilde{m}} \bigcap cQ_{0,m} \neq \emptyset}}\mathcal{E}_{l}(f, cQ_{k,\widetilde{m}})_{r}, \quad k \in \mathbb{N}_{0}, m \in \mathbb{Z}^{n},
\end{equation}

\begin{equation}
\label{eq2.13}
\frac{1}{C}\delta^{l}_{r}(Q,Q)f \le \mathcal{E}_{l}(f,Q)_{r} \le C \delta^{l}_{r}(Q,Q)f.
\end{equation}

The inequality \eqref{eq2.12} in the case $r \geq 1$ follows from Theorem~2 of~\cite{Br}. In the general case the proof follows the same lines (with the use of
Markov's inequality for polynomials). Estimate \eqref{eq2.13} is a~very deep result; the first proof of it was published in~\cite{Br}
(the case $r \geq 1$) and in~\cite{St}  (the general setting).
Note that the constant $C  > 0$ in inequalities \eqref{eq2.12}, \eqref{eq2.13} depends only on $n,c,l,r$ but not on~$f$.

\begin{Def} \rm
\label{Def2.5}
Let $p,q,r,\sigma_{1},\sigma_{2} \in (0,\infty]$, $\alpha_{1},\alpha_{2} \in \mathbb{R}$, $\alpha_{3} \geq 0$,
and let $\{t_{k}\} \in \mathcal X^{\alpha_{3}}_{\alpha,\sigma,p}$ be a weight sequence.
By $\widetilde{B}^{l}_{p,q,r}(\{t_{k}\})$ we denote the linear space of all functions $f \in L^{\mathrm{loc}}_{r}$ with finite quasi-norm
\begin{equation}
\label{eq2.14}
\|f|\widetilde{B}^{l}_{p,q,r}(\{t_{k}\})\|:= \Bigl(\sum\limits_{k=0}^{\infty} \|t_{k} \delta^{l}_{r}(\cdot+2^{-k}I) f|L_{p}\|^{q}\Bigr)^{\frac{1}{q}}+\|t_{0}\|f|L_{r}(\cdot+I)\||L_{p}\|.
\end{equation}
(the modifications of \eqref{eq2.14} in the case $q=\infty$ are trivial).
\end{Def}

\begin{Remark} \rm
\label{Remark2.4}
Even thought one may formally define the space $\widetilde{B}^{l}_{p,q,r}(\{t_{k}\})$ for all $\alpha_{1},\alpha_{2} \in \mathbb{R}$,
it is however of possible interest only when $\alpha_{2} \geq \alpha_{1} \geq 0$. This space is complete for nonnegative $\alpha_{1},\alpha_{2}$. Besides,
for $l > \alpha_{2}$ and $\sigma_{2} \geq p$ it is independent of~$l$, the corresponding norms being equivalent. For proofs, see~\cite{Ty1} (Sections~2 and~4).
\end{Remark}

\begin{Remark} \rm
\label{Remark2.5}
According to \cite{Ty1},
\begin{gather*}
\|f|\widetilde{B}^{l}_{p,q,r}(\{t_{k}\})\| \approx \|f|\widetilde{B}^{l}_{p,q,r}(\{\overline{t}_{k}\})\| \approx \\
\approx \Bigl(\sum\limits_{k=0}^{\infty}\Bigl(\sum\limits_{m \in \mathbb{Z}^{n}}t^{p}_{k,m}\mathcal{E}_{l}(f,cQ_{k,m})_{r}\Bigr)^{\frac{q}{p}}\Bigr)^{\frac{1}{q}}+\Bigl(\sum\limits_{m \in \mathbb{Z}^{n}}t^{p}_{0,m}\|f|L_{r}(Q_{0,m})\|^{p}\Bigr)^{\frac{1}{p}},
\end{gather*}
where $c \geq 1$ is an arbitrary parameter. Here, one norm can be estimated in terms of the other one with a~constant depending on the parameter~$c$ and on~$\alpha_{3},n,l,r,p,q$.
Hence, using~\eqref{eq2.8} it easily follows that, for any $k_{0} \in \mathbb{N}_{0}$,
\begin{equation}
\label{eq2.15}
\|f|\widetilde{B}^{l}_{p,q,r}(\{t_{k}\})\| \approx \Bigl(\sum\limits_{k=k_{0}}^{\infty}\Bigl(\sum\limits_{m \in \mathbb{Z}^{n}}t^{p}_{k,m}\mathcal{E}_{l}(f,cQ_{k,m})_{r}\Bigr)^{\frac{q}{p}}\Bigr)^{\frac{1}{q}}+\Bigl(\sum\limits_{m \in \mathbb{Z}^{n}}t^{p}_{k_{0},m}\|f|L_{r}(Q_{k_{0},m})\|^{p}\Bigr)^{\frac{1}{p}}.
\end{equation}
In \eqref {eq2.15} the constant through which one norm is estimated in terms of the other depends only on $c, n, l, k_{0}, \alpha, \alpha_{3}, \sigma, r,p,q$.
\end{Remark}

\begin{Def} \rm
\label{Def2.6}
Let $U$ be an open subset of $\mathbb{R}^{n}$, $p,q,r,\sigma_{1},\sigma_{2} \in (0,\infty]$, $\alpha_{1},\alpha_{2} \in \mathbb{R}$, $\alpha_{3} \geq 0$
and let  $\{t_{k}\} \in \mathcal X^{\alpha_{3}}_{\alpha,\sigma,p}$ be a~weight sequence. The spaces $\widetilde{B}^{l}_{p,q,r}(U,\{t_{k}\})$ and $B^{\varphi_{0}}_{p,q}(U,\{t_{k}\})$ are defined to be the restrictions of the corresponding spaces from $\mathbb{R}^{n}$ to $U$. They are endowed with the quotient space quasi-norms. More precisely, for $f \in L^{\mathrm{loc}}_{r}(U)$,
$$
\|f|\widetilde{B}^{l}_{p,q,r}(U,\{t_{k}\})\| = \inf\{\|g|\widetilde{B}^{l}_{p,q,r}(\{t_{k}\})\|:g|_{U} = f \text{ a.e. on   } U\},
$$
and, for $f \in D'(U)$,
$$
\|f|B^{\varphi_{0}}_{p,q}(U,\{t_{k}\})\|=\inf\{\|g|B^{\varphi_{0}}_{p,q}(\{t_{k}\})\|:g|_{U} = f \text{ in the sense of } D'(U)\}.
$$
\end{Def}

\begin{Remark} \rm
\label{Remark2.6}
Let $p \in (0,\infty)$,  $\alpha_{1}=\alpha_{2}=s \in \mathbb{R}$, and let $\gamma^{p} \in A^{\mathrm{loc}}_{\infty}$ be a~weight.
Putting $\{t_{k}\}=\{2^{ks}\gamma\}$ in Definitions  \ref{Def2.4}, \ref{Def2.5}, we get at the definitions of the weighted Besov spaces (see~\cite{Ry}),
which will be denoted by $\widetilde{B}^{s}_{p,q,r}(\gamma)$, $B^{s}_{p,q}(\gamma)$, respectively. Indeed, we note that if for $\nu \in (1,\infty)$
a~weight $\gamma^{p} \in A^{\textrm{loc}}_{\nu}$, then for  $s \in \mathbb{R}$ the weight sequence $\{t_{k}\}=\{2^{ks}\gamma\}$ lies in
$\mathcal{X}^{\alpha_{3}}_{\alpha,\sigma,p}$
with suitable parameters $\alpha_{3}\alpha,\sigma$ (see \cite{Ty1}, \cite{Ty2} for details). We deliberately drop the index~$l$
in the notation of the weighted Besov spaces $\widetilde{B}^{s}_{p,q,r}(\gamma)$, because from the results of~\cite{Ty1} it follows that for
$\gamma^{p} \in A^{\mathrm{loc}}_{\frac{p}{r}}$ and $l > s > 0$ the corresponding spaces coincide, the norms being equivalent.
We also skip the symbol $\varphi_{0}$ in the notation for the spaces $B^{s}_{p,q}(\gamma)$, because the norms corresponding to different functions $\varphi_{0}$ are equivalent
for $\gamma \in A^{\mathrm{loc}}_{\infty}$ and under suitable  conditions on the parameter $L_{\varphi}$  (see \cite{Ry} for details).
\end{Remark}

\textbf{Problem A.}
Find an intrinsic description of the spaces $B^{\varphi_{0}}_{p,q}(U,\{t_{k}\})$ and $\widetilde{B}^{l}_{p,q,r}(U,\{t_{k}\})$. In
other words, it is required to find equivalent norms in the  spaces $B^{\varphi_{0}}_{p,q}(U,\{t_{k}\})$ and $\widetilde{B}^{l}_{p,q,r}(U,\{t_{k}\})$,
which would utilize only the information about the distribution (function) on an open set~$U$.

\begin{Remark} \rm
\label{Remark2.7}
Problem~A  is closely related with the problem of constructing a~bounded extension operator from the spaces $B^{\varphi_{0}}_{p,q}(U,\{t_{k}\})$ and
$\widetilde{B}^{l}_{p,q,r}(U,\{t_{k}\})$ into the spaces $B^{\varphi_{0}}_{p,q}(\{t_{k}\})$ and $\widetilde{B}^{l}_{p,q,r}(\{t_{k}\})$, respectively.
Below we shall show that as an extension operator for the space $B^{\varphi_{0}}_{p,q}(U,\{t_{k}\})$
one may use the Rychkov operator~\cite{Ry2}, which is linear.
To extend functions from the space $\widetilde{B}^{l}_{p,q,r}(U,\{t_{k}\})$ we shall need the nonlinear operator that was employed in~\cite{DeVore}
in solving a~similar problem for classical Besov spaces.
\end{Remark}

\begin{Remark} \rm
\label{Remark2.8}
In view of Remark~\ref{Remark2.4} the solution to Problem~{A} for spaces of variable smoothness leads automatically to the
solution of the corresponding problems for weighted Besov spaces with Muckenhoupt weights, which are new results.
\end{Remark}

We shall require below the following theorem. This is the classical Hardy inequality for sequences. For a~proof see, for example,~\cite{Bennet}.

\begin{Th}
\label{Th2.1}
Let $0 < q \le \infty$, $0 < \mu \le q$, $\alpha \geq 0$, and let $\{a_{k}\}$ and $\{b_{k}\}$ be two sequences of real numbers.
Then the inequality (with constant $C > 0$ independent of sequences $\{a_{k}\}$, $\{b_{k}\}$)
$$
\Bigl(\sum\limits_{k=0}^{\infty}2^{k\alpha q}|b_{k}|^{q}\Bigr)^{\frac{1}{q}} \le C \Bigl(\sum\limits_{k=0}^{\infty}2^{k\alpha q}|a_{k}|^{q}\Bigr)^{\frac{1}{q}}
$$
holds under the following conditions:

{\rm 1)} $\lambda > \alpha$ and
$$
|b_{k}| \le C 2^{-k\lambda}\Bigl(\sum\limits_{j=0}^{k}2^{j\lambda}|a_{j}|^{\mu}\Bigr)^{\frac{1}{\mu}},
$$

{\rm 2)} \hfill $\displaystyle
|b_{k}| \le C \Bigl(\sum\limits_{j=k+1}^{\infty}|a_{j}|^{\mu}\Bigr)^{\frac{1}{\mu}}$.\hfill \null
\end{Th}

\section{Traces of the spaces $B^{\varphi_{0}}_{p,q}(\{t_{k}\})$ on Lipschitz domains}

Throughout this section we fix  parameters $p,q, \sigma_{1}, \sigma_{2} \in (0,\infty]$, $\alpha_{1},\alpha_{2} \in \mathbb{R}$, $\alpha_{3} \geq 0$
and a $p$-ad\-mis\-si\-ble weight sequence $\{t_{k}\} \in \mathcal X^{\alpha_{3}}_{\alpha,\sigma,p}$.

Like in \cite{Ry} for a~Lipschitz domain, we shall consider either a special or a~bounded Lipschitz domain.

A special Lipschitz domain is defined as an open set $G$ lying above the graph of a Lipschitz function. More precisely,
$$
G=\{x=(x',x_{n}) \in \mathbb{R}^{n}: x_{n} > \varpi(x')\},
$$
where $\varpi$ is a function satisfying the Lipschitz condition on $\mathbb{R}^{n-1}$ with constant $M \geq 1$. We set $K:=\{x=(x',x_{n})|x_{n} > M|x'|\}$.

A bounded Lipschitz domain is a bounded domain $G$, whose boundary $\partial G$ can be covered by a~finite number of balls $B_{k}$ so that, possibly after a~proper rotation,
$\partial G \bigcap B_{k}$ for each $k$ is a~part of the graph of a~Lipschitz function.

\begin{Lm}
\label{Lm3.1}
Let $\omega \in D$, $f \in B^{\varphi_{0}}_{p,q}(\{t_{k}\})$. Next, let $\varphi_{0} \in D(-K)$ have nonzero integral and let $\varphi(x)=\varphi_{0}(x)-2^{-n}\varphi_{0}(\frac{x}{2})$, $1+L_{\varphi} > \alpha_{2}$, $\sigma_{2} \geq p$. Then $\omega f \in B^{\varphi_{0}}_{p,q}(\{t_{k}\})$ and
\begin{equation}
\label{eq3.1}
\|\omega f|B^{\varphi_{0}}_{p,q}(\{t_{k}\})\| \le C \|f|B^{\varphi_{0}}_{p,q}(\{t_{k}\})\|,
\end{equation}
where the constant $C > 0$ in independent of~$f$.
\end{Lm}

\begin{Remark}
\label{Rem3.1}\rm
Lemma \ref{Lm3.1} enables one to easily reduce the problem on the intrinsic description of Besov spaces of variable smoothness on a~bounded Lipschitz domain~$G$
to a~similar problem on a~special Lipschitz domain~$G$. It will be convenient to give the proof of Lemma~\ref{Lm3.1} at the end of this section.
\end{Remark}

In view of Remark~\ref{Rem3.1}, we shall assume throughout this section that $G$ is a~fixed special Lipschitz domain.

\begin{Lm} (the local reproducing formula)
\label{Lm3.2}
Let $G$ be a special Lipschitz domain, let $\varphi_{0} \in D(-K)$ have nonzero integral, and let $\varphi(x)=\varphi_{0}(x)-2^{-n}\varphi_{0}(\frac{x}{2})$. Then,
for any given $L \in \mathbb{N}_{0}$, there exist functions $\psi_{0},\psi \in D(-K)$ such that $L_{\psi} \geq L$ and
\begin{equation}
\label{eq3.2}
 f = \sum\limits_{k=0}^{\infty}\psi_{k}\ast\varphi_{k}\ast f \quad  \text{ in } D'(G)
\end{equation}
for all $f \in D'(G)$.
\end{Lm}

For a proof, see \cite{Ry2}.

\begin{Remark} \rm
\label{Rem3.2}
Of  course, the local reproducing formula similar  to \eqref{eq3.2} also holds for distributions  $f \in D'$ (see \cite{Ry}, \cite{Ry2} for details).
\end{Remark}

Let $\{g^{k}\}$ be a sequence of measurable functions. Given  $j \in \mathbb{N}_{0}$, $A > 0$, $m \in \mathbb{Z}^{n}$, we set
$$
G^{j}_{A}(m):=\sup\limits_{k \geq j}2^{A(j-k)}\sup\limits_{y \in Q_{j,m}}|g^{k}(y)|.
$$

Let $A > 0$,  a~function $\varphi_{0} \in D$, $\int\varphi_{0}(x)\,dx=1$, $\varphi=\varphi_{0}-2^{-n}\varphi_{0}(\frac{\cdot}{2})$ and $c > 1$.
For a~distribution $f \in S'_{e}$, we consider the maximal function
\begin{multline}
\label{eq3.3}
M_{A,\varphi_{0}}(m,j,c)[f]:=M_{A}(m,j,c)[f]:=\sup\limits_{k \geq j}2^{A(j-k)} \sup\limits_{y \in cQ_{j,m}} |\varphi_{k} \ast f (y)| \\
\hbox{ for}\
m \in \mathbb{Z}^{n}, \ j \in \mathbb{N}_{0},\  c \geq 1.
\end{multline}

Let us denote by $S'_{e}(G)$ the subset of $D'(\Omega)$ consisting of distributions having finite order and at most exponential growth
at infinity, that is, $f \in S'_{e}(G)$ if and only if the estimate

$$
|<f,\theta>| \le C(f) \sup\limits_{x \in G, |\alpha| \le N(f)}|D^{\alpha}f|\hbox{exp}(N(f)|x|), \quad \hbox{ for all } \theta \in D(G).
$$

For a distribution  $f \in S'_{e}(G)$, we  consider the maximal function
\begin{multline}
\label{eq3.4}
M^{G}_{A,\varphi_{0}}(m,j,c)[f]:=M^{G}_{A}(m,j,c)[f]:=\sup\limits_{k \geq j}2^{A(j-k)} \sup\limits_{y \in cQ_{j,m} \bigcap G} |\varphi_{k} \ast f (y)| \\
 \hbox{ for } m \in \mathbb{Z}^{n}, \ j \in \mathbb{N}_{0}, \ c \geq 1.
\end{multline}

The maximal functions $M_{A,\varphi_{0}}(m,j,c)[f]$, $M^{G}_{A,\varphi_{0}}(m,j,c)[f]$ provide a~convenient tool in replacing the Peetre maximal functions.
However, it is worth noting that Rychkov~\cite{Ry} introduced functions which are close to our maximal functions $M_{A,\varphi_{0}}(m,j,c)[f]$, $M^{G}_{A,\varphi_{0}}(m,j,c)[f]$
(but which differ in possessing no localization). The advantage of our maximal functions is that they are local and prove to be very convenient  (unlike Peetre-type maximal functions)
in obtaining estimates in the spaces $B^{\varphi_{0}}_{p,q}(\{t_{k}\})$. Roughly speaking, they enable one to get rid of `exponentially' decreasing tails---the fact which
substantially simplifies the technicalities.

\begin{Remark} \rm
\label{Rem3.3}
We note that $M^{G}_{A}(m,j,c) < \infty$, provided that $f \in S'_{e}(G)$. The proof of this fact follows the same lines as that of Lemma 2.9 in~\cite{Ry}.
\end{Remark}

The next lemma plays a crucial role in the proof of some results that follow.

\begin{Lm}
\label{Lm3.3}
Let $c \geq 1$, $\varphi_{0} \in D$, $\int \varphi_{0}(x) \,dx$ and $\varphi:=\varphi_{0}-\frac{1}{2^{n}}\varphi_{0}(\frac{\cdot}{2})$. Assume that $A > \max\{-\alpha_{1},0\}$. Then, for any distribution $f \in S'_{e}(G)$,
\begin{equation}
\begin{split}
\label{eq3.5}
   \Bigl(\sum\limits_{j=0}^{\infty}\Bigl(\sum\limits_{m \in \mathbb{Z}^{n}}t^{p}_{j,m} (M^{G}_{A}(m,j,c)[f])^{p}\Bigr)^{\frac{q}{p}}\Bigr)^{\frac{1}{q}} \le C\Bigl(\sum\limits_{k=0}^{\infty}\Bigl(\int\limits_{G}t^{p}_{k}(x)(\varphi_{k} \ast f(x))^{p}\,dx \Bigr)^{\frac{q}{p}}\Bigr)^{\frac{1}{q}},
\end{split}
\end{equation}
(the modifications in the case $p=\infty$ or $q=\infty$ are straightforward).
Here, the constant $C:=C(n,p,q,\alpha_{3},\alpha_{1},\sigma_{1},c,A,\varphi_{0}) > 0$ is independent of  $f$.
\end{Lm}

The proof is similar to that of Lemma 3.3 in~\cite{Ty2}, but for the sake of completeness we give the details.
We consider the case $p,q \neq \infty$ (the case $p=\infty$ or $q=\infty$ is dealt with similarly).

We claim that if $A,r > 0$, then, for $j \in \mathbb{N}_{0}, m \in \mathbb{Z}^{n}$, $c \geq 1$,
\begin{equation}
\label{eq3.6}
(M^{G}_{A}(m,j,c)[f])^{r} \le C \sum\limits_{k=j}^{\infty}2^{(j-k)Ar}2^{kn}\int\limits_{cQ_{j,m} \bigcap G }|\varphi_{k} \ast f(z)|^{r}\,dz.
\end{equation}
Here, the constant $C$ on the right of \eqref{eq3.6} is independent both of $x,j,m$ and the function~$f$. The derivation of \eqref{eq3.6} depends on
a~variant of Stromberg--Torchinsky's trick, which was used in the proof of Lemma~2.9 in~\cite{Ry}.

In view of \eqref{eq3.2}, we have
$$
\varphi_{j}\ast f = (\psi_{0})_{j} \ast (\varphi_{0})_{j} \ast \varphi_{j}\ast f + \sum\limits_{k=j+1}^{\infty} \varphi_{j} \ast \psi_{k} \ast \varphi_{k} \ast f,
$$
where one may assume that $L_{\psi} > A$.

For all $k \geq j$ we have the following estimate (which depends on the condition  $L_{\psi} \geq A$, see~\cite{Ry} for more details)
\begin{equation}
\label{eq3.7}
\|\varphi_{j} \ast \psi_{k} | L_{\infty}\| \le C 2^{(j-k)A+jn},
\end{equation}
in which the constant  $C > 0$ depends on $c,L_{\psi},\varphi_{0}$, but is independent of both $k$ and~$j$.

Hence, since for $k \geq j$ the function $\varphi_{j} \ast \psi_{k}$  has support inside a~cube of side length at most $\widetilde{c} 2^{-jn}$,
we have, for $i \le j$,
\begin{equation}
\label{eq3.8}
\sup\limits_{y \in c_{1}Q_{i,m} \bigcap G} \varphi_{j}\ast f(y) \le C \Bigl(\sum\limits_{k=j}^{\infty}2^{(j-k)A}2^{kn}\int\limits_{c_{2}Q_{i,m}\bigcap G}|\varphi_{k}\ast f(z)|\,dz\Bigr).
\end{equation}

In the case $r \geq 1$, the proof concludes by applying the H\"older inequality first for the integrals and then for series with exponents $r,r'$.

In the case $r \in (0,1)$, it clearly follows from\eqref{eq3.8} that, for $k \geq j$,
\begin{equation}
\begin{split}
\label{eq3.9}
&2^{(j-k)A}\sup\limits_{y \in c_{1}Q_{j,m} \bigcap G} \varphi_{k}\ast f(y) \le C \Bigl(\sum\limits_{l=k}^{\infty}2^{(j-l)A}2^{ln}\int\limits_{c_{2}Q_{j,m}\bigcap G}|\varphi_{l}\ast f(z)|\,dz\Bigr) \le \\
&\le C (M^{G}_{A}(m,j,c_{1}))^{1-r}\Bigl(\sum\limits_{l=j}^{\infty}2^{(j-l)Ar}2^{ln}\int\limits_{c_{2}Q_{j,m}\bigcap G}|\varphi_{l}\ast f(z)|^{r}\,dz\Bigr).
\end{split}
\end{equation}

Estimate \eqref{eq3.6} follows from \eqref{eq3.9}, provided that  $M^{G}_{A}(m,j,c_{1})[f] < \infty$.

To show that $M^{G}_{A}(m,j,c_{1})[f]$ is finite, we note that, for a distribution $f \in S'_{e}(G)$,
$$
|\varphi_{k} \ast f(x)| \le C\hbox{exp}(N_{f}|x|)p_{N_{f}}(\varphi_{k}).
$$

It follows that $M^{G}_{A}(m,j,c_{1})[f] < \infty$ for $A > C(N_{f})$, and hence estimate \eqref{eq3.6} holds with  $A > C(N_{f})$.
As a~result, since the right-hand side of~\eqref{eq3.6} increases with decreasing~$A$, we have
\begin{equation}
\label{eq3.10}
\sup\limits_{y \in cQ_{j,m} \bigcap G}|\varphi_{j} \ast f(y)| \le C(N_{f})\Bigl(\sum\limits_{l=j}^{\infty}2^{(j-l)Ar}2^{ln}\int\limits_{c_{2}Q^{n}_{j,m}\bigcap G}|\varphi_{l}\ast f(z)|^{r}\,dz\Bigr)^{\frac{1}{r}}
\end{equation}
for $A,r > 0$ (now for any values!) and under the condition that the right-hand side of~\eqref{eq3.6} is finite.

The constant $C > 0$ on the right of \eqref{eq3.10}  depends on $N_{f}$, because the corresponding constant on the right of \eqref{eq3.7} depends on the
parameter $L_{\psi}$, which in turn depends on~$A$ (recall that we assumed $L_{\psi} \geq A > C(N_{f})$).

However, substituting \eqref{eq3.10} in the definition of $M^{G}_{A}(m,j,c)[f]$, we easily see that
\begin{equation}
\label{eq3.11}
 M^{G}_{A}(m,j,c)[f] \le C(N_{f})\Bigl(\sum\limits_{l=j}^{\infty}2^{(j-l)Ar}2^{ln}\int\limits_{c_{2}Q^{n}_{j,m} \bigcap G}|\varphi_{l}\ast f(z)|^{r}\,dz\Bigr)^{\frac{1}{r}} < \infty,
\end{equation}
provided that the right-hand side in \eqref{eq3.6} is finite.

Now estimate \eqref{eq3.11} in combination with \eqref{eq3.8}, \eqref{eq3.9} (in which the constant $C > 0$ is independent of $f$, $j$,~$m$) proves
\eqref{eq3.6} with the constant $C > 0$ independent of $f$, $j$,~$m$.

It is clear that we can choose $r \in (0,\infty]$ such that $\sigma_{1}=r p'_{r}$. We now employ estimate \eqref{eq3.6}, and then use the monotonicity of $l_{q}$ in  $q$, apply the Minkowski inequality for sums (because
$\frac{p}{\mu} \geq 1$),  the
 H\"older inequality for integrals with exponents $p_{r}:=\frac{p}{r}$, $p'_{r}$, and finally use condition \eqref{eq2.6} with $\sigma_{1}=r p'_{r}=\frac{pr}{p-r}$.
As a~result, we have, for $j \in \mathbb{N}_{0}$ and $\mu \le \min\{1,q,r\}$,
\begin{gather}
\Bigl(\sum\limits_{m \in \mathbb{Z}^{n}}t^{p}_{j,m} [M^{G}_{A}(m,j,c)]^{p}\Bigr)^{\frac{\mu}{p}} \le  \\
\le
C \Bigl(\sum\limits_{m \in \mathbb{Z}^{n}}t^{p}_{j,m} \Bigl(\sum\limits_{k=j}^{\infty}2^{(j-k)Ar}2^{kn}\int\limits_{\widetilde{c} Q_{j,m} \bigcap G}|\varphi_{k}\ast f(z)|^{r}\,dz\Bigr)^{\frac{\mu p}{\mu r}}\Bigr)^{\frac{ r}{p}} \le
\notag \\
\le C \Bigl(\sum\limits_{m \in \mathbb{Z}^{n}}t^{p}_{j,m} \Bigl(\sum\limits_{k=j}^{\infty}\Bigl(2^{(j-k)Ar+kn}\int\limits_{\widetilde{c} Q_{j,m}\bigcap G}|\varphi_{k}\ast f(z)|^{r}\,dz\Bigr)^{\frac{\mu}{r}}\Bigr)^{\frac{p}{\mu}}\Bigr)^{\frac{\mu}{p}} \le
\notag \\
\le C \sum\limits_{k=j}^{\infty}2^{(j-k)A\mu}\Bigl(2^{\frac{knp}{r}}\sum\limits_{m \in \mathbb{Z}^{n}}t^{p}_{j,m}\Bigl(\int\limits_{\widetilde{c} Q_{j,m} \bigcap G}\frac{t^{r}_{k}(x)}{t^{r}_{k}(x)}|\varphi_{k}\ast f(z)|^{r}\,dz\Bigr)^{\frac{p}{r}} \Bigr)^{\frac{\mu}{p}} \le
\notag \\
\label{eq3.12}
\le C \sum\limits_{k=j}^{\infty}2^{(j-k)(A+\alpha_{1}) \mu}\Bigl(\sum\limits_{m \in \mathbb{Z}^{n}}\int\limits_{\widetilde{c}Q_{j,m} \bigcap G}t_{k}^{p}(z)|\varphi_{k}\ast f(z)|^{p}\,dz\Bigr)^{\frac{\mu}{p}}.
\end{gather}

Here, we also used the fact that $2^{\frac{knp}{r}}=2^{kn}+2^{kn\frac{p}{\sigma_{1}}}$.

Now the required assertion follows from \eqref{eq3.12} and Theorem~\ref{Th2.1} with the above conditions on the parameters $A$ and $\alpha_{1}$ .

\begin{Lm}
\label{Lm3.4}
The space $B^{\varphi_{0}}_{p,q}(\{t_{k}\})$ is a~complete quasi-normed vector space which
is continuously embedded into the space $S'_{e}$.
\end{Lm}

\textbf{Proof}. We first show that $B^{\varphi_{0}}_{p,q}(\{t_{k}\}) \subset S'_{e}$ and  that the embedding operator is continuous.
By the definition of the space $S'_{e}$, it will suffice for our purposes to establish the estimate
\begin{equation}
\begin{split}
\label{eq3.13}
|\langle f,\theta\rangle| \quad \le \quad C \|f|B^{\varphi_{0}}_{p,q}(\{t_{k}\})\|\sup\limits_{x \in \mathbb{R}^{n}}\Bigl(\sum\limits_{|\alpha| \le L}|D^{\alpha}\theta(x)|\hbox{exp}(N|x|)\Bigr)
\end{split}
\end{equation}
with $f \in B^{\varphi_{0}}_{p,q}(\{t_{k}\})$ and $\theta \in D$, in which
the constants $C, N, L$ depend on $\{t_{k}\}$, $p$,~$q$,~$n$, but are independent of both  $f$ and~$\theta$.

Estimate \eqref{eq3.13} is easily seen to follow from the estimate
\begin{equation}
\begin{split}
\label{eq3.14}
|\widetilde{\theta} \ast f(x)|  \le C \|f|B^{\varphi_{0}}_{p,q}(\{t_{k}\})\|\Bigl(\sup\limits_{\widetilde{x} \in \mathbb{R}^{n}}\sum\limits_{|\alpha| \le L}|D^{\alpha}\widetilde{\theta}(\widetilde{x})|\Bigr)\hbox{exp}(N|x|),
\end{split}
\end{equation}
in which $\widetilde{\theta} \in D(I)$, the constant $C$ is independent of both  $f$ and $\widetilde{\theta}$.

To prove \eqref{eq3.14} we shall argue as in the derivation of estimate \eqref{eq3.6},
replacing the function $\varphi_{0}$ by the function $\widetilde{\theta}$. We also note that the constant~$C$ on the right of \eqref{eq3.7} (with $\varphi_{j}$ replaced by  $\widetilde{\theta}$)
can easily be estimated from above by the number $\sup\limits_{\widetilde{x} \in \mathbb{R}^{n}}\sum\limits_{|\alpha| \le L} |D^{\alpha}\widetilde{\theta}(\widetilde{x})|$
(which follows from Taylor's formula with remainder in the Lagrange form). Now using Remark \ref{Remark2.2}, Lemma \ref{Lm3.3} and the estimate
$\sup\limits_{x \in Q_{0,m}}|\theta \ast f(x)| \le M_{A}(m,0,c)$ with $m \in \mathbb{Z}^{n}$, $A > \max\{-\alpha_{1},0\}$, we obtain (note that $L$ depends on $A$)
\begin{gather}
\sup\limits_{x \in Q_{0,m}}|\widetilde{\theta} \ast f(x)| \le \frac{1}{t_{0,m}}t_{0,m}\sup\limits_{x \in Q_{0,m}}|\widetilde{\theta} \ast f(x)| \le \notag \\
\le C  \hbox{exp}(\widehat{N}(\alpha_{3})|m|)\Bigl(\sum\limits_{m \in \mathbb{Z}^{n}}t^{p}_{0,m}\sup\limits_{x \in Q_{0,m}}|\widetilde{\theta} \ast f(x)|^{p}\Bigr)^{\frac{1}{p}} \le
\notag \\
\label{eq3.15}
 \le C \|f|B^{\varphi_{0}}_{p,q}(\{t_{k}\})\| \Bigl(\sup\limits_{x \in \mathbb{R}^{n}}\sum\limits_{|\alpha| \le L}|D^{\alpha}\widetilde{\theta}(x)|\Bigr) \hbox{exp}(\widehat{N}(\alpha_{3})|m|).
\end{gather}

As a result, from \eqref{eq3.15} we get \eqref{eq3.14} with $N=N(\widehat{N}(\alpha_{3}),n)$.
The proof of the completion of the space $B^{\varphi_{0}}_{p,q}(\{t_{k}\})$ repeats verbatim the arguments from the last paragraph of Lemma~2.15 in \cite{Ry},
estimate \eqref{eq3.13} being useful.

This proves the lemma.

\begin{Remark}
\label{Rem3.4} \rm
Assume that $\{t_{k}\} \in \mathcal{W}^{\alpha_{3}}_{\alpha_{1},\alpha_{2}}$ for some $\alpha_{3} \geq 0$, $\alpha_{1},\alpha_{2} \in \mathbb{R}$. Then we have
$$
t_{0,m} \geq C (1+|m|)^{-d}
$$

for some $C > 0$ and $d$ depending on $\alpha_{3}$. The argument used to deduce \eqref{eq3.14} in our case in fact gives \eqref{eq3.13} with $(1+|x|)^{d}$ instead of $\hbox{exp}(N|x|)$. This shows that if $f \in S'_{e}$ belongs to $B^{\varphi_{0}}_{p,q}(\{t_{k}\})$ with $\{t_{k}\} \in \mathcal{W}^{\alpha_{3}}_{\alpha_{1},\alpha_{2}}$ then $f$ belongs to $B^{\varphi_{0}}_{p,q}(\{t_{k}\}) \bigcap S'$ and hence $f \in B^{\{t_{k}\}}_{p,q}$ (here we used Definition \ref{Def2.2} theorem 3.10 from \cite{Kempka4} and \ref{Remark2.4}).
\end{Remark}

\begin{Lm}
\label{Lm3.5}
Let  $\{g^{j}\}$ be a sequence of measurable functions. Let $\varphi_{0} \in D$, $\int \varphi_{0}(x)\,dx =1$, $\varphi=\varphi_{0}-2^{-n}\varphi_{0}(\frac{\cdot}{2})$ and $L_{\varphi}+1 > \alpha_{2}$, $\sigma_{2} \geq p$, $A > 0$.
Next, let  $\psi_{0}$ be the same as in~\eqref{eq3.2} and $L_{\psi}+1 > A$. Assume that
$$
\sum\limits_{j=0}^{\infty}\Bigl(\sum\limits_{m \in \mathbb{Z}^{n}} t^{p}_{j,m}\Bigl(G^{j}_{A}(m)\Bigr)^{p}\Bigr)^{\frac{q}{p}} < \infty
$$
(the modifications in the case $p=\infty$ or $q=\infty$ are clear).

Then the series $\sum\limits_{j=0}^{\infty} \psi_{j}\ast g^{j}$ converges in $S'_{e}$. Moreover,
\begin{equation}
\label{eq3.16}
\|\sum\limits_{j=0}^{\infty}\psi_{j}\ast g^{j}|B^{\varphi_{0}}_{p,q}(\{t_{k}\})\| \le C \Bigl(\sum \limits_{j=0}^{\infty}\Bigl(\sum\limits_{m \in \mathbb{Z}^{n}} t^{p}_{j,m}\Bigl(G^{j}_{A}(m)\Bigr)^{p}\Bigr)^{\frac{q}{p}}\Bigr)^{\frac{1}{q}}
\end{equation}
(the modifications in the case $p=\infty$ or $q=\infty$ are evident).
\end{Lm}

\textbf{Proof}. In view of Remark 2.1 we may assume without loss of generality that $\sigma_{2}=p$.

We claim that $\psi_{j} \ast g^{j} \in S'_{e}$, provided the right-hand side of  \eqref{eq3.16} is finite.

Indeed, using Remark 2.2 for $\omega \in S_{e}$ we have
$$
|\langle \psi_{j} \ast g^{j}, \omega \rangle | = \left|\int \omega(x) \psi_{j} \ast g^{j}(x) \,dx\right| \le \sum\limits_{m \in \mathbb{Z}^{n}}\int\limits_{Q_{j,m}}|\omega(x) \psi_{j} \ast g^{j}(x)| \,dx \le
$$

$$
\le 2^{-jn}\sup\limits_{m \in \mathbb{Z}^{n}}\frac{t_{j,m}}{t_{j,m}}\sup\limits_{x \in Q_{j,m}}|\omega(x) \psi_{j} \ast g^{j}(x)| \le $$

$$\le C 2^{-jn} \left(\sup\limits_{x \in \mathbb{R}^{n}}|\omega(x)| \text{exp} (N(j,\alpha_{3})|x|)\right)\left(\sup\limits_{m \in \mathbb{Z}^{n}}t_{j,m}G^{j}_{A}(m)\right) \le C p_{N(j,\alpha_{3})}[\omega].
$$

In \cite{Ry} it was shown that
\begin{equation}
\label{eq3.17}
\|\varphi_{l} \ast \psi_{j}|L_{\infty}\| \le C
\begin{cases}
 2^{(l-j)(L_{\psi}+1)}2^{ln}, & l \le j \\
 2^{(j-l)(L_{\varphi}+1)}2^{jn}, & j < l.
\end{cases}
\end{equation}

The diameter of the support of the function $\varphi_{l} \ast \psi_{j}$ is at most  $c 2^{-jn}$ in the case $l > j$ and is at most
$c 2^{-ln}$ in the case $j \geq l$ (here, the constant~$c$ depends only on $\varphi_{0}$, $\psi_{0}$, $n$). Combining this fact with estimate \eqref{eq3.17}, this establishes
\begin{equation}
\label{eq3.18}
\sup\limits_{y \in Q_{l,m}}|\varphi_{l} \ast \psi_{j} \ast g^{j}(y)| \le C
\begin{cases}
 2^{(l-j)(L_{\psi}+1)}\sup\limits_{y \in cQ_{l,m}}|g^{j}(y)|, & l \le j \\
 2^{(j-l)(L_{\varphi}+1)}\sup\limits_{y \in cQ_{j,\widetilde{m}}}|g^{j}(y)|, & j < l.
\end{cases}
\end{equation}

On the right of \eqref{eq3.18} the cube $Q_{j,\widetilde{m}}$,  $j < l$, is the only dyadic cube of side length $2^{-j}$ that contains the cube $Q_{l,m}$

For $l > j$, using \eqref{eq3.18} and (2.5), we find that
\begin{equation}
\begin{split}
\label{eq3.19}
&\int t^{p}_{l}(x)|\varphi_{l} \ast \psi_{j} \ast g^{j}(x)|^{p}\,dx \le C 2^{p(j-l)(L_{\varphi}+1)}\sum\limits_{\widetilde{m} \in \mathbb{Z}^{n}}\sum\limits_{\substack{m \in \mathbb{Z}^{n} \\ Q_{l,m} \subset Q_{j,\widetilde{m}}}}t^{p}_{l,m} \sup\limits_{y \in cQ_{j,\widetilde{m}}}|g^{j}(y)|^{p} \le \\
&\le C 2^{p(j-l)(L_{\varphi}+1-\alpha_{2})} \sum\limits_{\widetilde{m} \in \mathbb{Z}^{n}}t^{p}_{j,\widetilde{m}} \sup\limits_{y \in cQ_{j,\widetilde{m}}}|g^{j}(y)|^{p} \le C 2^{p(j-l)(L_{\varphi}+1-\alpha_{2})} \sum\limits_{m \in \mathbb{Z}^{n}}t^{p}_{j,m}\Bigl(G^{j}_{A}(m)\Bigr)^{p}.
\end{split}
\end{equation}

If $l \le j$, then by \eqref{eq3.18} we have
\begin{equation}
\begin{split}
\label{eq3.20}
&\int t^{p}_{l}(x)|\varphi_{l} \ast \psi_{j} \ast g^{j}(x)|^{p}\,dx \le C 2^{p(l-j)(L_{\psi}+1)}\sum\limits_{m \in \mathbb{Z}^{n}}t^{p}_{l,m}\sup\limits_{y \in cQ_{l,m}}|g^{j}(y)|^{p} \le \\ &\le C2^{p(l-j)(L_{\psi}+1-A)}\sum\limits_{m \in \mathbb{Z}^{n}}t^{p}_{l,m}\Bigl(G^{l}_{A}(m)\Bigr)^{p}
\end{split}
\end{equation}

Assuming that   $\min\{L_{\psi}+1-A, L_{\varphi}+1-\alpha_{2}\} > 0$, we set  $\varepsilon:=\max\{L_{\psi}+1-A, L_{\varphi}+1-\alpha_{2}\}$. Then
$$
\|\psi_{j} \ast g^{j}| B^{\varphi_{0}}_{p,q}(\{2^{-2k\varepsilon}t_{k}\})\|^{q} \le 2^{-j\varepsilon}\Bigl(\sum\limits_{m \in \mathbb{Z}^{n}} t^{p}_{j,m}\Bigl(G^{j}_{A}(m)\Bigr)^{p}\Bigr)^{\frac{q}{p}}.
$$

Hence, by Lemma \ref{Lm3.4} the series $\sum\limits_{j=0}^{\infty}\psi_{j} \ast g^{j}$ converges in the complete
space $B^{\varphi_{0}}_{p,q}(\{2^{-2 k\sigma}t_{k}\})$, and hence, in $S'_{e}$ (again by Lemma \ref{Lm3.4}) to some distribution~$g$.

If $\min\{L_{\psi}+1-A, L_{\varphi}+1-\alpha_{2}\} > 0$, then from \eqref{eq3.19}, \eqref{eq3.20} with
$\mu \in (0,\min\{1,p,q\}]$ and $l \in \mathbb{N}_{0}$ we get  (since the $l_{q}$-norm is monotone in~$q$)
\begin{gather}
\Bigl(\int t^{p}_{l}(x) |\varphi_{l} \ast g(x)|\Bigr)^{\frac{1}{p}} \le \notag\\
\le C \Bigl(\sum\limits_{m \in \mathbb{Z}^{n}}t^{p}_{l,m}\sup\limits_{y \in Q_{l,m}}|\varphi_{l} \ast g(y)|^{p}\Bigr)^{\frac{1}{p}} \le \notag \\
C \Bigl(\sum\limits_{m \in \mathbb{Z}^{n}}t^{p}_{l,m}\Bigl(\sum\limits_{j=0}^{l}2^{\mu(j-l)(L_{\psi}+1)}\sup\limits_{y \in cQ_{j,\widetilde{m}}}|g^{j}(y)|^{\mu}\Bigr)^{\frac{p}{\mu}}\Bigr)^{\frac{1}{p}}+
\notag\\
\label{eq3.21}
+C \Bigl(\sum\limits_{m \in \mathbb{Z}^{n}}t^{p}_{l,m}[G^{l}_{A}(m)]^{p}\Bigr)^{\frac{1}{p}}=:S_{1,l}+S_{2,l}.
\end{gather}

To estimate $S_{1,l}$ we shall employ the Minkowski inequality for sums (since $\frac{p}{\mu} \geq 1$) and condition~\eqref{eq2.5}. We have
\begin{gather}
S_{1,l} \le C \Bigl(\sum\limits_{m \in \mathbb{Z}^{n}}t^{p}_{l,m}\Bigl(\sum\limits_{j=0}^{l}2^{\mu(j-l)(L_{\varphi}+1)}\sup\limits_{y \in cQ_{j,m}}|g^{j}(y)|^{\mu}\Bigr)^{\frac{p}{\mu}}\Bigr)^{\frac{1}{p}} \le
\notag\\
\le C \Bigl(\sum\limits_{j=0}^{l}2^{\mu(j-l)(L_{\varphi}+1)}\Bigl(\Bigl(\sum\limits_{\substack{\widetilde{m} \in \mathbb{Z}^{n}\\
Q_{l,m} \subset Q_{j,\widetilde{m}}}} t^{p}_{l,m}\Bigr)\sup\limits_{y \in cQ_{j,\widetilde{m}}}|g^{j}(y)|^{p}\Bigr)^{\frac{\mu}{p}}\Bigr)^{\frac{1}{\mu}} \le
\notag\\
\le C \Bigl(\sum\limits_{j=0}^{l}2^{\mu(l-j)(\alpha_{2}-L_{\varphi}-1)}\Bigl(\sum\limits_{\widetilde{m} \in \mathbb{Z}^{n}}t^{p}_{j,\widetilde{m}}\sup\limits_{y \in cQ_{j,\widetilde{m}}}|g^{j}(y)|^{p}\Bigr)^{\frac{\mu}{p}}\Bigr)^{\frac{1}{\mu}} \le
\notag\\
\label{eq3.22}
\le C \Bigl(\sum\limits_{j=0}^{l}2^{\mu(l-j)(\alpha_{2}-L_{\varphi}-1)}\Bigl(\sum\limits_{m \in \mathbb{Z}^{n}}t^{p}_{j,m}\Bigl(G^{j}_{A}(m)\Bigr)^{p}\Bigr)^{\frac{\mu}{p}}\Bigr)^{\frac{1}{\mu}}.
\end{gather}

The conclusion of the lemma now follows from Theorem~\ref{Th2.1}, estimates \eqref{eq3.21}, \eqref{eq3.22}, and the constraint  $L_{\varphi}+1>\alpha_{2}$.
\smallskip

Now we are ready to prove Lemma \ref{Lm3.1}.

\textbf{Proof of Lemma 3.1}. We give only a~sketch proof, because
many of the details are similar to those from~\cite{Ry} (the proof of Theorem~2.21).

First of all, in view of Remark~\ref{Remark2.3} we may without loss of generality assume that the function $\varphi_{0}$ is
chosen so that $\varphi$ has the required (sufficiently big) number of zero moments (this number will be fixed at the end of the proof).

Assume first that $f \in L_{1}^{\mathrm{loc}}$. Applying
Taylor's formula with integral remainder to the function~$\omega$, this gives
\begin{equation}
\begin{split}
\label{eq3.23}
&\varphi_{j} \ast (\omega f)(x)= \sum\limits_{|\beta| \le T-1}\frac{(-1)^{|\beta|}}{\beta!}D^{\beta}\omega(x)\int y^{\beta}\varphi_{j}(y)f(x-y)\,dy+\\
&+\sum\limits_{|\beta|=T}\frac{(-1)^{T}}{\beta!}\int y^{\beta} \varphi_{j}(y)\Bigl[\int\limits_{0}^{1}D^{\beta}\omega(x-\tau y)(1-\tau)^{T-1}\,d\tau\Bigr]f(x-y)\,dy.
\end{split}
\end{equation}

Let now $f \in S'_{e}$. We set $f^{l}:=\psi_{l} \ast \varphi_{l} \ast f$ for $l \in \mathbb{N}_{0}$. Using the local reproducing formula (Lemma \ref{Lm3.2} and Remark \ref{Rem3.2}),
we see that $\sum\limits_{l=0}^{k}f^{l} \to f$ in $D'$ as $k \to \infty$. Hence, taking into account that
 $\sum\limits_{l=0}^{k}f^{l} \in C^{\infty} \subset L^{\mathrm{loc}}_{1}$, from
 \eqref{eq3.23} and since the operator of multiplication by a~smooth function in~$D'$ is continuous, we find that
\begin{equation}
\begin{split}
\label{eq3.24}
&\varphi_{j} \ast (\omega f)(x) = \sum\limits_{l=1}^{\infty} \sum\limits_{|\beta| \le T-1}\frac{(-1)^{|\beta|}}{\beta!}D^{\beta}\omega(x)\int y^{\beta}\varphi_{j}(y)f^{l}(x-y)\,dy+\\
&+\sum\limits_{l=1}^{\infty}\sum\limits_{|\beta|=T}\frac{(-1)^{T}}{\beta!}\int y^{\beta} \varphi_{j}(y)\left[\int\limits_{0}^{1}D^{\beta}\omega(x-\tau y)(1-\tau)^{T-1}\,d\tau\right]f^{l}(x-y)\,dy =\\
&= \Sigma_{j,1}(x)+\Sigma_{j,2}(x), \quad x \in \mathbb{R}^{n}.
\end{split}
\end{equation}

Making an obvious change of variables, it follows from Fubini's theorem that
\begin{equation}
\begin{split}
\label{eq3.25}
|\Sigma_{j,2}(x)| \le C \sum\limits_{l=0}^{\infty}\sum\limits_{|\beta|=T}\int \left|\int y^{\beta} \varphi_{j}(y)J_{\beta}(x,y)\psi_{l}(z-y)\varphi_{l} \ast f (x-z)\,dy\right|dz,
\end{split}
\end{equation}
where $J_{\alpha}(x,y)$ is the expression in the square brackets in~\eqref{eq3.24}.

We have two cases to consider: $l \le j$ and $l > j$.

If $l \le j$, then using the inclusions $\operatorname{supp}\psi_{l} \subset 2^{-l}I$, $\operatorname{supp}\varphi_{j} \subset 2^{-j}I$ and
the estimate $\sup\limits_{y}|J_{\beta}(x,y)\psi_{l}(z-y)| \le C 2^{ln}$, we find that
\begin{multline}
\label{eq3.26}
\sum \limits_{|\beta|=T}\int \int \left|y^{\beta} \varphi_{j}(y)J_{\beta}(x,y)\psi_{l}(z-y)\right| |\varphi_{l} \ast f (x-z)|\,dy\,dz \le \\
\le C 2^{-jT} \sup\limits_{x' \in x+c2^{-l}I}|\varphi_{l} \ast f(x')|
\end{multline}

 If  $l > j$, then the absolute value of the inner integral on the right of \eqref{eq3.25} is estimated from above by
 $C 2^{L_{\psi}(j-l)}2^{jn}\sup\limits_{x' \in x+c2^{-j}I}|\varphi_{l} \ast f(x')|$ (for details, see \cite{Ry}, the derivation of estimate (2.50)). Hence,
\begin{multline}
\sum\limits_{|\beta|=T}\int \left|\int y^{\beta} \varphi_{j}(y)J_{\beta}(x,y)\psi_{l}(z-y) \varphi_{l} \ast f (x-z)\,dy\right|dz \le
\\
\label{eq3.27}
C 2^{-jT}2^{L_{\psi}(j-l)} \sup\limits_{x' \in x+c2^{-j}I}|\varphi_{l} \ast f(x')|
\end{multline}

Employing \eqref{eq3.26} and \eqref{eq3.27}, this gives
\begin{equation}
\begin{split}
\label{eq3.28}
&t^{p}_{j,m}\sup\limits_{x \in Q_{j,m}}|\Sigma_{j,2}(x)|^{p} \le C t^{p}_{j,m}\Bigl(\sum\limits_{l=0}^{j} 2^{-jT}\sup\limits_{x \in c Q_{l,\widetilde{m}}}|\varphi_{l} \ast f(x)|\Bigr)^{p}+\\
&+C t^{p}_{j,m}2^{-jpT}\Bigl(\sum\limits_{l=j}^{\infty} 2^{(L_{\psi}-A)(j-l)}\Bigr)^{p}\Bigl(M_{A}(m,j,c)[f]\Bigr)^{p},
\end{split}
\end{equation}
provided that $L_{\psi} > A$

Using \eqref{eq3.28} and arguing as in the derivation of estimate \eqref{eq3.21}, we obtain
\begin{multline}
\label{eq3.29}
\Bigl(\sum\limits_{m \in \mathbb{Z}^{n}}t^{p}_{j,m}\sup\limits_{x \in Q_{j,m}}|\Sigma_{j,2}(x)|^{p}\Bigr)^{\frac{1}{p}} \le
\\
\le C  \Bigl(\sum\limits_{l=0}^{j}2^{\mu(l-j)T}\Bigl(\sum\limits_{\widetilde{m} \in \mathbb{Z}^{n}}t^{p}_{l,\widetilde{m}} \sup\limits_{x \in c Q_{l,\widetilde{m}}}|\varphi_{l} \ast f(x)|^{p}\Bigr)^{\frac{\mu}{p}}\Bigr)^{\frac{1}{\mu}}+\\
+C 2^{-jT} \Bigl(\sum\limits_{m \in \mathbb{Z}^{n}}t^{p}_{j,m}\Bigl(M_{A}(j,m,c)[f]\Bigr)^{p}\Bigr)^{\frac{1}{p}}.
\end{multline}

Next, by a~similar argument as in the derivation of \eqref{eq3.22}, it follows from \eqref{eq3.29}, Lemma \ref{Lm3.3} and Theorem \ref{Th2.1} that
\begin{equation}
\begin{split}
\label{eq3.30}
\sum\limits_{j=0}^{\infty}\Bigl(\sum\limits_{m \in \mathbb{Z}^{n}}t^{p}_{j,m}\sup\limits_{x \in Q_{j,m}}|\Sigma_{j,2}(x)|^{p}\Bigr)^{\frac{1}{p}} \le C \|f|B^{\varphi_{0}}_{p,q}(\{t_{k}\})\|,
\end{split}
\end{equation}
provided that $T > \alpha_{2}$.

To estimate $\Sigma_{1,j}(x)$ we note that
\begin{equation}
\label{eq3.31}
|\Sigma_{1,j}(x)| \le C \sum\limits_{|\beta| \le T-1}2^{-j|\beta|}|\zeta^{\beta}_{j} \ast f(x)|, \quad x \in \mathbb{R}^{n}.
\end{equation}
where $\zeta^{\beta}(y)=y^{\beta}\varphi(y)$ with  $y \in \mathbb{R}^{n}$.

By estimate (3.13) of~\cite{Ty2}, we have
\begin{equation}
\begin{split}
\label{eq3.32}
&\Bigl(\sum\limits_{m \in \mathbb{Z}^{n}}t^{p}_{j,m}\sup\limits_{x \in Q_{j,m}}|\Sigma_{j,1}(x)|^{p}\Bigr)^{\frac{1}{p}} \le\\
& \le C \sum\limits_{|\beta| \le T-1} \Bigl(\sum\limits_{l=0}^{j}2^{\mu(l-j)(L_{\zeta^{\beta}}+1-|\beta|)}\Bigl(\sum\limits_{\widetilde{m} \in \mathbb{Z}^{n}}t^{p}_{l,\widetilde{m}} \sup\limits_{x \in c Q_{l,\widetilde{m}}}|\varphi_{l} \ast f(x)|^{p}\Bigr)^{\frac{\mu}{p}}\Bigr)^{\frac{1}{\mu}}+\\
&+C \sum\limits_{|\beta| \le T-1}2^{-j|\beta|} \Bigl(\sum\limits_{m \in \mathbb{Z}^{n}}t^{p}_{j,m}\Bigl(M_{A}(j,m,c)[f]\Bigr)^{p}\Bigr)^{\frac{1}{p}}.
\end{split}
\end{equation}

If the function $\varphi$ has zero moments up to the order~$L$, then the function $\zeta^{\beta}$ has zero moments up to the order $L-|\beta|$.
If $L_{\varphi} > T+\alpha_{2}$, then similarly to~\eqref{eq3.30} we get
\begin{equation}
\begin{split}
\label{eq3.33}
\sum\limits_{j=0}^{\infty}\Bigl(\sum\limits_{m \in \mathbb{Z}^{n}}t^{p}_{j,m}\sup\limits_{x \in Q_{j,m}}|\Sigma_{j,1}(x)|^{p}\Bigr)^{\frac{1}{p}} \le C \|f|B^{\varphi_{0}}_{p,q}(\{t_{k}\})\|.
\end{split}
\end{equation}

Estimates \eqref{eq3.30} and \eqref{eq3.33}  conclude the proof of the lemma.

\begin{Remark}
\label{Remark3.4} \rm
Lemma \ref{Lm3.1} extends the right-hand side of Theorem~4.3 in \cite{Kempka5} (in the case of constant $p,q$) to the case of more general weight sequences $\{t_{k}\}$.
\end{Remark}

\begin{Th}
\label{Th3.1}
Let $\varphi_{0} \in D(-K)$, $\int \varphi_{0}(x)\,dx = 1$, $\varphi:=\varphi_{0}-2^{-n}\varphi_{0}(\frac{\cdot}{2})$. Let $\psi_{0}, \psi$ be the same as in \eqref{eq3.2} and $1+L_{\psi} > \max{0,-\alpha_{1}}$. Let $\sigma_{2} \geq p$, $L_{\varphi}+1 > \alpha_{2}$.  Then the map $\text{Ext}: D'(G) \mapsto D'$, as defined by
\begin{equation}
\label{eq3.34}
\text{Ext}[f]=\sum\limits_{j=0}^{\infty}\psi_{j} \ast (\varphi_{j} \ast f)_{G}
\end{equation}
defines a~linear bounded operator from the space $B^{\varphi_{0}}_{p,q}(G,\{t_{k}\})$ into the space $B^{\varphi_{0}}_{p,q}(\{t_{k}\})$.

Furthermore,
\begin{equation}
\label{eq3.35}
\|f|B^{\varphi_{0}}_{p,q}(G,\{t_{k}\})\| \approx \Bigl(\sum\limits_{j=0}^{\infty}\Bigl(\sum\limits_{m \in \mathbb{Z}^{n}} t^{p}_{j,m}\Bigl(M^{G}_{A}(m,j,c)[f]\Bigr)^{p}\Bigr)^{\frac{q}{p}}\Bigr)^{\frac{1}{q}}
\end{equation}
(the modifications in the case $p=\infty$ or $q=\infty$ are straightforward).
\end{Th}

\textbf{Proof}. If $f \in B^{\varphi_{0}}_{p,q}(G,\{t_{k}\})$ then $g \in S'_{e}$ (see Definition \ref{Def2.6}). Clearly, for every $g \in S'_{e}$ such that $g=f$ on $G$ (in the sence of $D'(G)$) and for $c \geq 1$, $A > 0$, we have
$$
M^{G}_{A}(m,j,c)[f] \le M_{A}(m,j,c)[g], \quad j \in \mathbb{N}_{0}, m \in \mathbb{Z}^{n}.
$$

Hence, using Lemma \ref{Lm3.3}(with $\mathbb{R}^{n}$ instead of $G$)  and Definition \ref{Def2.6},
$$
\Bigl(\sum\limits_{j=0}^{\infty}\Bigl(\sum\limits_{m \in \mathbb{Z}^{n}} t^{p}_{j,m}\Bigl(M^{G}_{A}(m,j,c)[f]\Bigr)^{p}\Bigr)^{\frac{q}{p}}\Bigr)^{\frac{1}{q}} \le C \|f|B^{\varphi_{0}}_{p,q}(G,\{t_{k}\})\|.
$$

To conclude the proof of Theorem~\ref{Th3.1}  it suffices to establish the opposite estimate.

To this end we note that, for $\{g^{j}\}=\{(\varphi_{j} \ast f)_{G}\}$,
$$
G^{j}_{A}(m)=M^{G}_{A}(m,j,1)[f], \quad j \in \mathbb{N}_{0}, m \in \mathbb{Z}^{n}
$$

Hence, using Lemma \ref{Lm3.4} and Definition \ref{Def2.6},
$$
\|f|B^{\varphi_{0}}_{p,q}(G,\{t_{k}\})\| \le \|\operatorname{Ext}[f]|B^{\varphi_{0}}_{p,q}(\{t_{k}\})\| \le C \Bigl(\sum\limits_{j=0}^{\infty}\Bigl(\sum\limits_{m \in \mathbb{Z}^{n}} t^{p}_{j,m}\Bigl(M^{G}_{A}(m,j,c)[f]\Bigr)^{p}\Bigr)^{\frac{q}{p}}\Bigr)^{\frac{1}{q}}.
$$

This completes the proof of  Theorem~\ref{Th3.1}.

\begin{Remark}
\label{Remark3.5} \rm
Theorem \ref{Th3.1} is an intermediate step towards the solution of Problem~\textbf{A} for the spaces $B^{\varphi_{0}}_{p,q}(\{t_{k}\})$. Of course,
Theorem~\ref{Th3.1} gives an intrinsic description of the trace. The necessary and sufficient conditions on a~trace are expressed in terms
which depend on the information about the \_\_distribution~$f$\_\_ only on the domain~$G$.
At the same time, the verification of the test in Theorem~3.1 is a~challenge from the numerical point of view.
\end{Remark}

Our next theorem gives a pretty simple description of the trace space.

\begin{Th}
\label{Th3.2}
Let $\varphi_{0} \in D(-K)$, $\int \varphi_{0}(x)\,dx = 1$ and $\varphi:=\varphi_{0}-2^{-n}\varphi_{0}(\frac{\cdot}{2})$. Let $\sigma_{2} \geq p$ and $L_{\varphi}+1 > \alpha_{2}$. Then, for all $f \in S'_{e}(G)$,
\begin{equation}
\label{eq3.36}
\|f|B^{\varphi_{0}}_{p,q}(G,\{t_{k}\})\| \approx \Bigl(\sum\limits_{j=0}^{\infty}\Bigl(\int\limits_{G} t^{p}_{j}(x)|\varphi_{j} \ast f(x)|^{p}\,dx\Bigr)^{\frac{q}{p}}\Bigr)^{\frac{1}{q}}.
\end{equation}
\end{Th}

\textbf{Proof} The estimate `$\geq$' clearly follows from the definition of the norm in the space $B^{\varphi_{0}}_{p,q}(G,\{t_{k}\})$.
The reverse estimate $\le$ is secured by Theorem~3.1 and Lemma~\ref{Lm3.3}.

\section{Traces of the spaces $\widetilde{B}^{l}_{p,q,r}(\{t_{k}\})$ on $(\varepsilon,\delta)$-domains}

The Lipschitz domains (considered in the previous section)
constitute a~very special case of $(\varepsilon, \delta)$-domains. It is worth noting that $(\varepsilon, \delta)$-domains may have
highly irregular boundary  (see~\cite{J}). The fundamental paper of Jones~\cite{J} was succeeded by papers pertaining to
the problem of extension of functions from various function spaces with $(\varepsilon, \delta)$-domains on the entire
space $\mathbb{R}^{n}$ \cite{Chua}, \cite{DeVore}, \cite{Vodop}, \cite{Chr}, \cite{Se},\cite{Mi},\cite{Mi2}. As we shall see below, the machinery
developed in~\cite{DeVore} for the study of extensions of classical Besov spaces can be applied, after certain modifications, to work out the theory of Besov spaces of variable smoothness on
$(\varepsilon, \delta)$-domains.

\begin{Def}\rm  (\cite{J})
\label{Def4.1}
Let $\varepsilon, \delta > 0$. An open set $\Omega \subset \mathbb{R}^{n}$ is called an $(\varepsilon,\delta)$-\textit{domain} if, for any $x,y \in  \Omega$ such that $|x-y| < \delta$, there exists
a~rectifiable path $\Gamma \subset \Omega$ of length $< \frac{1}{\varepsilon}|x-y|$ connecting $x$ and $y$ such that, for each $z \in \Gamma$,
$$
\operatorname{dist}(z,\partial \Omega) > \varepsilon \frac{|z-x| |z-y|}{|x-y|}.
$$
\end{Def}

We fix throughout this section some $(\varepsilon,\delta)$-domain $\Omega$. We shall also assume that  $\operatorname{rad}(\Omega)=\inf
\limits_{\alpha}\inf\limits_{x \in \Omega_{\alpha}}\sup\limits_{y \in \Omega_{\alpha}}|x-y| > 0$. Here,
$\Omega_{\alpha}$ denotes the connected components of an open set~$\Omega$.

We shall frequently use the following notations. For any cube~$Q$, we set $Q^{\ast}:=\frac{9}{8}Q$.
By $F$ and $F_{c}$ we shall denote, respectively, the Whitney decomposition of the open sets
$\Omega$ and $\mathbb{R}^{n} \setminus \overline{\Omega}$ (see~\cite{DeVore} for details).

Following \cite{DeVore}, for a~cube $Q \in F_{c}$, we let $Q^{s}$ denote any cube from $F$ of~maximal diameter such that
$\operatorname{dist}(Q^{s},Q) < 2 \operatorname{dist}(Q, \partial \Omega)$. The cube $Q^{s}$ will be called the reflection of $Q$. By~$\mathcal{F}_{c}$
we shall denote the cubes from the family $F_{c}$ whose diameters are at most~$\delta$.

\begin{Lm}
\label{Lm4.1}
For every cube $Q \in \mathcal{F}_{c}$ we have

{\rm 1)} \begin{equation}
\label{eq4.1}
 \frac{1}{C(n,\varepsilon,\delta)}l(Q) \le l(Q^{s}) \le C(n,\varepsilon,\delta) l(Q);
   \end{equation}

{\rm 2)} \begin{equation}
   \label{eq4.2}
   \operatorname{dist}(Q,Q^{s})  \le C(n,\varepsilon,\delta) l(Q);
   \end{equation}

{\rm 3)} each cube $Q^{s} \in F$ is a reflected cube for at most $C(n,\varepsilon,\delta)$ cubes from the family $\mathcal{F}_{c}$.
\end{Lm}

The proof of this lemma is contained in \S\,5 of~\cite{DeVore} (see the arguments succeeding estimate (5.3)).

Following \cite{DeVore}, we construct a~nonlinear (since the local near-best approximations depend on~$f$) extensions of the operator from
the $(\varepsilon,\delta)$-domain $\Omega$ to the entire $\mathbb{R}^{n}$.
Let $r \in (0,1]$ and assume that $f \in L_{r}(Q \bigcap \Omega)$ for every cube $Q$. We set
\begin{equation}
\label{eq4.3}
\widetilde{f}:=\operatorname{Ext}[f]=\chi_{\Omega}f+\sum\limits_{\substack{Q \in F_{c} \\ \hbox{diam}(Q) \le \delta}}P_{Q^{s}}[f]\phi_{Q},
\end{equation}
where $\{\phi_{Q}\}_{Q \in F_{c}}$ is the partition of unity for the open set $\mathbb{R}^{n} \setminus \overline{\Omega}$,
and $P_{Q^{s}}[f]$ is a~polynomial of near-best approximation with constant $\lambda \geq 1$ (see formula (2.8) in \cite{DeVore} for details) to a~function~$f$  in the quasi-norm $L_{r}(Q^{s})$ (see \cite{DeVore} for more details).

Actually,  \eqref{eq4.3} defines  a  family  of extension  operators,  since  each  choice  of near
best  approximants  $P_{Q^{s}}[f]$  give  an  extension  $\widetilde{f}$.  The  results  that  follow  apply  to
any such  extension operator $\operatorname{Ext}$ with the restriction  that  the constant $\lambda \geq 1$ is  fixed

The proof of the main result depends on two auxiliary lemmas.

\begin{Lm} (\cite{DeVore})
\label{Lm4.2}
Let $R_{0}$ and $Q$ be two cubes from $F$ with $\hbox{diam}(Q) \le \hbox{diam}(R_{0})$
and $\operatorname{dist}(Q, R_{0}) \le \min \{\delta, C_{1}\hbox{diam}(R_{0})\}$ with $C_{1}$ a fixed constant. Then, there
is a sequence of cubes $Q =: R_{m}, R_{m-1}, \dotsc , R_{0}$, from $F$, such that each $R_{j}$
touches $R_{j-1}$, $j = l, \dotsc, m$, and for each $j = l, \dotsc, m$, $R_{j} \subset c R_{0}$ and for
each $j = 0, . .. , m - 1$, $Q \subset c R_{j}$ with $c > 0$ depending only on $C_{1}$ and $\Omega$.
\end{Lm}

\begin{Lm} (\cite{DeVore})
\label{Lm4.3}
Let $\beta > 0$ and $\operatorname{Ext}[f]$ be defined by \eqref{eq4.1}. Let $R$ be a cube with $\operatorname{dist}(R,\partial \Omega) \le \hbox{diam} R \le a\delta$ where $a$ is a fixed constant depending only on $\varepsilon,\delta$ and $n$. Let $f \in L_{r}(Q \bigcap \Omega)$ for every cube $Q$, $\beta \le r \le 1$. Then  we have
\begin{equation}
\label{eq4.4}
\left(E_{l}(\widetilde{f},R)_{r}\right)^{r} \le C \sum\limits_{\substack{S \in  F \\ S \subset cR}}\left(E_{l}(f,S^{\ast})_{r}\right)^{r}
\end{equation}
where $c,C$ depend only on $l,n,\beta,\lambda$ and $\varepsilon, \delta$.
\end{Lm}

For convenience, for $0<r \le p \le \infty$ we set $p_{r}=\frac{p}{r}$. The exponent $p'_{r}$ is determined from the equation $\frac{1}{p_{r}}+\frac{1}{p'_{r}}=1$.

The main result of this section is as follows

\begin{Th}
Let $0<\alpha_{1} \le \alpha_{2} < l$, $r \in (0,p]$, $0 < p,q \le \infty$, $p \geq r$, $\sigma_{1} \geq r p'_{r}$, $\sigma_{2} \geq p$.
Next, let $p$~be an admissible weight sequence $\{t_{k}\}$ such that  $\{\overline{t}_{k}\} \in \mathcal{X}^{\alpha_{3}}_{\alpha,\sigma,p}$. Then
\begin{equation}
\label{eq4.5}
\begin{split}
&\|f|\widetilde{B}^{l}_{p,q,r}(\Omega,\{t_{k}\})\| \approx \Bigl(\sum\limits_{k=1}^{\infty}\Bigl(\int\limits_{\Omega}t^{p}_{k}(x)\left[\delta^{l}_{r}(x+2^{-k}I,\Omega)f \right]^{p}\,dx\Bigr)^{\frac{q}{p}}\Bigr)^{\frac{1}{q}}+\\
&+\Bigl(\int\limits_{\Omega}t^{p}_{0}(x)\left\|f|L_{r}\Bigl((x+I)\bigcap \Omega\Bigr)\right\|^{p}\,dx\Bigr)^{\frac{1}{p}}.
\end{split}
\end{equation}
(the modifications in the cases $p=\infty$ or $q=\infty$ are clear.)
\end{Th}

\textbf{Proof}. The estimate `$ \geq $' in  \eqref{eq4.5} is clear. To establish the estimate `$\le$', we only consider the case $p,q \in (0,\infty)$,
because in the remaining cases the proof follows the same ideas, but is technically simpler.

For further purposes we note that in view of Remark~\ref{Remark2.1} we may assume without loss of generality that  $\sigma_{1}=r p'_{r}$, $\sigma_{2}=p$. Moreover in view of Corollary 4.4 \cite{Ty1} we have $\|\cdot|\widetilde{B}^{l}_{p,q,r}\|\approx \|\cdot|\widetilde{B}^{l}_{p,q,\tilde{r}}\|$ for $r < \tilde{r} \le p$ and $\sigma_{1} > \tilde{r} p'_{\tilde{r}} > r p'_{r}$. Hence, without loss of generality we may assume that $r \in (0,\min{1,p}]$.

Until the end of the proof we shall fix the number $k_{0} = \min\{k \in \mathbb{N} | \sqrt{n}2^{-k} < a\delta\}$ (the constant~$a$ is the same as in Lemma \ref{Lm4.3}).

In view of Remark~\ref{Remark2.5} it suffices to prove the estimate
\begin{equation}
\label{eq4.6}
\begin{split}
&\Bigl(\sum\limits_{k=1}^{\infty}\Bigl(\sum\limits_{m \in \mathbb{Z}^{n}}t^{p}_{k,m}\Bigl(\mathcal{E}_{l}(\widetilde{f},\widetilde{c}Q_{k,m})_{r}\Bigr)^{p}\Bigr)^{\frac{q}{p}}\Bigr)^{\frac{1}{q}}+
\Bigl(\sum\limits_{m \in \mathbb{Z}^{n}}t^{p}_{0,m}\left\|\widetilde{f}|L_{r}(Q_{k_{0},m})\right\|^{p}\Bigr)^{\frac{1}{p}} \le\\
&\le C \Bigl(\sum\limits_{k=1}^{\infty}\Bigl(\sum\limits_{\substack{m \in \mathbb{Z}^{n} \\ Q^{\ast}_{k,m} \subset \Omega}}t^{p}_{k,m}\Bigl(\mathcal{E}_{l}(f,Q^{\ast}_{k,m})_{r}\Bigr)^{p}\Bigr)^{\frac{q}{p}}\Bigr)^{\frac{1}{q}}+\Bigl(\sum\limits_{\substack{m \in \mathbb{Z}^{n} \\ Q_{0,m} \subset \Omega }}t^{p}_{0,m}\left\|f|L_{r}\Bigl(Q_{0,m}\bigcap \Omega\Bigr)\right\|^{p}\Bigr)^{\frac{1}{p}}
\end{split}
\end{equation}
in which the constant $C > 0$ is independent of  function~$f$, and  $\widetilde{c} \in (1,\frac{9}{8})$.

We fix $k \geq k_{0}$ and estimate $\mathcal{E}_{l}(\widetilde{f},\widetilde{c}Q_{k,m})_{r}$ for different $m \in \mathbb{Z}^{n}$. Similarly to~\cite{DeVore},
we have three different cases to consider.

\textit{In the first case}, we consider the cubes such that $\operatorname{dist}(Q_{k,m},\partial \Omega) \le \sqrt{n}2^{-k} $ (the constant $a$ is the same as in Lemma \ref{Lm4.3}).
We let $A_{1}(k)$ denote the corresponding set of indexes~$m$.

Let $\mu \in (0,\min\{1,p,q\}]$ be fixed. Then, clearly by \eqref{eq4.3} and since the $l_{q}$-norm is monotone in~$q$, we have for $m \in A_{1}(k)$

$$
E_{l}(\widetilde{f}, \widetilde{c} Q_{k,m})_{r} \le \Bigl(\sum\limits_{j=k}^{\infty}\Bigl(\sum\limits_{\substack{\widetilde{m} \in \mathbb{Z}^{n} \\ Q_{j,\widetilde{m}} \subset c\widetilde{c}Q_{k,m}}}E_{l}(f,Q_{j,\widetilde{m}}^{\ast})_{r}\Bigr)^{\mu}\Bigr)^{\frac{1}{\mu}}.
$$

 In the case when a~cube $Q_{j,\widetilde{m}}$ does not lie in the family~$F$ it will be useful
to write $E_{l}(f,Q_{j,\widetilde{m}})_{r}=0$ to have the shorthand for some inequalities.
Hence, using the Minkowski inequality (since  $\frac{p}{\mu} \geq 1$), we find that
\begin{multline}
\label{eq4.7}
\Bigl(\sum\limits_{m \in A_{1}(k)}t^{p}_{k,m}\Bigl(\mathcal{E}_{l}(\widetilde{f}, \widetilde{c}Q_{k,m})_{r}\Bigr)^{p}\Bigr)^{\frac{\mu}{p}} \le
\\
\le C \sum\limits_{j=k}^{\infty}\Bigl(\sum\limits_{m \in A_{1}(k)}t^{p}_{k,m}2^{\frac{knp}{r}}\Bigl(\sum\limits_{\substack{\widetilde{m} \in \mathbb{Z}^{n} \\ Q_{j,\widetilde{m}} \subset c\widetilde{c}Q_{k,m} }}E_{l}(f,Q^{\ast}_{j,\widetilde{m}})_{r}\Bigr)^{p}\Bigr)^{\frac{\mu}{p}}.
\end{multline}

Let us estimate the inside sum in the outer brackets on the right of~\eqref{eq4.7}. To this end we shall use the monotonicity of the $l_{q}$-norm in~$q$
(taking into account that $0 < r \le 1$), apply the H\"older inequality  (since $\frac{p}{r} \geq 1$). We have
\begin{multline}
\label{eq4.8}
\sum\limits_{m \in A_{1}(k)}t^{p}_{k,m}2^{\frac{knp}{r}}\Bigl(\sum\limits_{\substack{\widetilde{m} \in \mathbb{Z}^{n}  \\ Q_{j,\widetilde{m}} \subset c\widetilde{c}Q_{k,m} }}E_{l}(f,Q^{\ast}_{j,\widetilde{m}})_{r}\Bigr)^{p} \le
\\
\le
\sum\limits_{m \in A_{1}(k)}t^{p}_{k,m}2^{\frac{knp}{r}}\Bigl(\sum\limits_{\substack{\widetilde{m} \in \mathbb{Z}^{n} \\ Q_{j,\widetilde{m}} \subset c\widetilde{c}Q_{k,m} }}\Bigl(E_{l}(f,Q^{\ast}_{j,\widetilde{m}})_{r}\Bigr)^{r}\Bigr)^{\frac{p}{r}} \le \\
\le \sum\limits_{m \in A_{1}(k)}t^{p}_{k,m}2^{\frac{(k-j)np}{r}}\Bigl(\sum\limits_{\substack{\widetilde{m} \in \mathbb{Z}^{n} \\ Q_{j,\widetilde{m}} \subset c\widetilde{c}Q_{k,m}}}\frac{1}{t^{\sigma_{1}}_{j,\widetilde{m}}}\Bigr)^{\frac{p}{\sigma_{1}}}\sum\limits_{\substack{\widetilde{m} \in \mathbb{Z}^{n} \\  Q_{j,\widetilde{m}} \subset c\widetilde{c}Q_{k,m} }}t^{p}_{j,\widetilde{m}}2^{\frac{jnp}{r}}(E_{l}(f,Q^{\ast}_{j,\widetilde{m}})_{r})^{p}.
\end{multline}
where the constant $c > 1$ is the same as in Lemma \ref{Lm4.1}.

From \eqref{eq4.7}, \eqref{eq4.8} in view of \eqref{eq2.6}, \eqref{eq2.7} (note that $\frac{n}{p}+\frac{n}{\sigma_{1}}=\frac{n}{r}$ and
that the overlap multiplicity of the cubes $c\widetilde{c}Q_{k,m}$ is finite with fixed~$k$ and variable~$m$) we have
\begin{equation}
\label{eq4.9}
\sum\limits_{m \in A_{1}(k)}t^{p}_{k,m}\Bigl(\mathcal{E}_{l}(\widetilde{f}, \widetilde{c}Q_{k,m})_{r}\Bigr)^{p} \le C \Bigl(\sum\limits_{j=k}^{\infty}2^{\mu(k-j)\alpha_{1}}\Bigl(\sum\limits_{\substack{\widetilde{m} \in \mathbb{Z}^{n} \\ Q^{\ast}_{j,\widetilde{m}} \subset \Omega }}t^{p}_{j,\widetilde{m}}\Bigl(\mathcal{E}_{l}(f,Q^{\ast}_{j,\widetilde{m}})_{r}\Bigr)^{p}\Bigr)^{\frac{\mu}{p}}\Bigr)^{\frac{p}{\mu}}
\end{equation}

In the \textit{second case}, we have $\operatorname{dist}(Q_{k,m},\partial \Omega) > \sqrt{n}2^{-k}$ and $\widetilde{c}Q_{k,m} \subset \Omega$.
So, $\widetilde{f}=\operatorname{Ext}[f] \equiv f$ and there is nothing to prove.

 In the \textit{third case}, $\operatorname{dist}(Q_{k,m},\partial \Omega) > \sqrt{n}2^{-k}$ and $Q_{k,m} \subset \Omega^{c}$.
 We let $A_{3}(k)$ denote the set of indexes~$m$ in this case. Clearly, for  $m \in A_{3}(k)$ the cube $Q_{k,m}$ is contained in a~unique Whitney cube $
 R \in F_{c}$. We have two cases to consider: in the first case $\operatorname{dist}(R,\partial \Omega) < a\delta$, and in the second case
 $\operatorname{dist}(R,\partial \Omega) \geq a\delta$. Accordingly, we split the index set $A_{3}(k)$ into two disjoint subsets
 $A^{1}_{3}(k)$ and $A^{2}_{3}(k)$. Clearly, $\sum\limits_{\substack{m \in \mathbb{Z}^{n}\\ \mbox{diam}(Q) < \delta}}\phi_{Q} = 1$ on~$R$
 if $\operatorname{dist}(R,\partial \Omega) < a\delta$ and $a$~is sufficiently small (the parameter~$a$ can
always be reduced without sacrificing the validity of Lemma~\ref{Lm4.3}). Hence, for $m \in A^{1}_{3}(k)$ and $Q_{k,m} \subset R$, $l(R) = 2^{-j}$, $j \in \mathbb{N}_{0}$,
we have
\begin{multline}
\label{eq4.10}
\Bigl(\delta^{l}_{r}(\widetilde{c}Q_{k,m},\widetilde{c}Q_{k,m})\Bigr)^{r} \le 2^{-jkr}\sup\limits_{|\alpha|=l}\sup\limits_{x \in \widetilde{c}Q_{k,m}}|D^{\alpha}\widetilde{f}(x)|^{r} \le
\\
\le C 2^{r(j-k)l}\sup \limits_{Q \in \Lambda_{R}}|P_{R^{s}}[f]-P_{Q^{s}}[f]|,
\end{multline}
where $\Lambda_{R}$ is the set of cubes from the family $F_{c}$ which have common boundary points with the cube~$R$.

By Lemma \ref{Lm4.1} and Lemma \ref{Lm4.2}, for any cube $Q \in \Lambda_{R}$ there exists a~chain of pairewisely touching cubes $R^{s}:=R_{0},\dotsc,Q^{s}$.
Besides, each cube from this chain lies in the family~$F$, and the quantity of cubes in this chain is bounded from above by $C(n,\varepsilon,\delta)$.
We let $T_{Q}$ denote this chain of cubes. An application of estimate 4.26 from~\cite{DeVore} shows that
\begin{equation}
\label{eq4.11}
|P_{R^{s}}[f]-P_{Q^{s}}[f]| \le C 2^{\frac{jn}{r}} \Bigl(\sum\limits_{S \in T_{Q}}E_{l}(f,S^{\ast})_{r}\Bigr)^{\frac{1}{r}}.
\end{equation}

Using \eqref{eq4.10}, \eqref{eq4.11}, this establishes
\begin{equation}
\label{eq4.12}
\mathcal{E}_{l}(\widetilde{f}, cQ_{k,m})^{r}_{r} \le C 2^{r(j-k)l}\sum\limits_{Q \in \Lambda_{R}} \sum\limits_{S \in T_{Q}}\mathcal{E}_{l}(f,S^{\ast})^{r}_{r}.
\end{equation}

We raise the left and right-hand side of \eqref{eq4.12} to the power $\frac{p}{r}$ and apply the
H\"older inequality to the right-hand side (because $\frac{p}{r} \geq 1$). We also note that
on the right of~\eqref{eq4.12} the number of terms is finite and is bounded above by a~constant, which depends only on $n,\varepsilon,\delta$. As a~result, we have
\begin{equation}
\label{eq4.13}
(\mathcal{E}_{l}(\widetilde{f}, cQ_{k,m})_{r})^{p} \le C 2^{p(j-k)l}\sum\limits_{Q \in \Lambda_{R}} \sum\limits_{S \in T_{Q}}\Bigl(\mathcal{E}_{l}(f,S^{\ast})_{r}\Bigr)^{p}.
\end{equation}

Next, using \eqref{eq4.13},
\begin{equation}
\label{eq4.14}
\sum\limits_{\substack{m \in A^{1}_{3}(k) \\ Q_{k,m} \subset R}}t^{p}_{k,m}\Bigl(\mathcal{E}_{l}(\widetilde{f}, Q_{k,m})_{r}\Bigr)^{p} \le C 2^{p(j-k)l}\Bigl(\sum\limits_{\substack{m \in A^{1}_{3}(k) \\ Q_{k,m} \subset R}}t^{p}_{k,m}\Bigr)\sum\limits_{Q \in \Lambda_{R}} \sum\limits_{S \in T_{Q}}\Bigl(\mathcal{E}_{l}(f,S^{\ast})_{r}\Bigr)^{p}.
\end{equation}

We now sum estimate \eqref{eq4.14} over all cubes $R \in \mathcal{F}_{c}$ that contain the cubes $Q_{k,m}$, $m \in A^{1}_{3}(k)$. In view of \eqref{eq2.5}, we have
\begin{equation}
\label{eq4.15}
\sum\limits_{m \in A^{1}_{3}(k)}t^{p}_{k,m} \mathcal{E}_{l}(\widetilde{f}, \widetilde{c}Q_{k,m})^{p}_{r} \le C \sum\limits_{S \in F}\Bigl(\frac{l(S)}{2^{k}}\Bigr)^{p(l-\alpha_{2})}t^{p}_{S}\Bigl(\mathcal{E}_{l}(f,S^{\ast})_{r}\Bigr)^{p}.
\end{equation}

If we now formally put $\mathcal{E}_{l}(f,Q_{j,m}^{\ast})_{r}=0$ for the cubes $Q_{j,m}$ not lying in the family~$F$,
 then estimate \eqref{eq4.15} can be rewritten in the form (recall that we fixed $k \geq k_{0}$)
\begin{multline}
\label{eq4.16}
\Bigl(\sum\limits_{m \in A^{1}_{3}(k)}t^{p}_{k,m} \Bigl(\mathcal{E}_{l}(\widetilde{f}, \widetilde{c}Q_{k,m})_{r}\Bigr)^{p}\Bigr)^{\frac{\mu}{p}\frac{1}{\mu}} \le
\\
\le C \Bigl( \sum\limits_{j=k_{0}}^{k}2^{\mu(l-\alpha_{2})(j-k)}\Bigl(\sum\limits_{\widetilde{m} \in \mathbb{Z}^{n}}t^{p}_{j,\widetilde{m}}\Bigl(\mathcal{E}_{l}(f,Q_{j,\widetilde{m}}^{\ast})_{r}\Bigr)^{p}\Bigr)^{\frac{\mu}{p}} \Bigr)^{\frac{1}{\mu}}.
\end{multline}
(here we also employ the condition $\mu \le p$ and use  the monotonicity of the $l_{q}$-norm in~$q$).

Arguing as in the proof of Theorem~5.3 of~\cite{DeVore}, one readily gets the following estimate for the
cubes $R \in \mathcal{F}_{c}$, $\operatorname{dist}(R,\partial \Omega) \geq a\delta$,
\begin{equation}
\label{eq4.17}
\sup\limits_{x \in R} |D^{\alpha}\widetilde{f}(x)| \le C \|f|L_{r}(R')\|^{p}
\end{equation}
where $R'$ is the union of all the cubes $Q^{s}$ such that $\phi_{Q}$ does not vanish on $R$.

Applying Taylor's formula with integral form of the remainder, it is immediate from \eqref{eq4.17} that
\begin{equation}
\label{eq4.18}
\Bigl(\delta^{l}_{r}(\widetilde{c}Q_{k,m}, \widetilde{c}Q_{k,m})\widetilde{f}\Bigr)^{p} \le C 2^{-klp}\|f|L_{r}(R')\|^{p}, \quad m \in A^{2}_{3}(k).
\end{equation}

Clearly, each set  $R'$ can be represented as a union of pairwise nonoverlapping  cubes $Q_{k_{0},m}$. Besides,
sets~$R'$ have finite overlap multiplicity (which depends only on $n,\varepsilon,\delta$), when $R$~runs over the cubes from~$\mathcal{F}_{c}$ for which
$\operatorname{dist}\{R, \partial \Omega\} \geq a\delta$. Using this fact, it follows from \eqref{eq2.5} and~\eqref{eq4.18} that
\begin{equation}
\label{eq4.19}
\sum\limits_{m \in A^{2}_{3}(k)}t^{p}_{k,m}\Bigl(\mathcal{E}_{l}(\widetilde{f}, \widetilde{c}Q_{k,m})_{r}\Bigr)^{p} \le C 2^{(k_{0}-k)lp}\sum\limits_{\substack{m \in \mathbb{Z}^{n} \\ Q_{k_{0},m} \subset \Omega}}t^{p}_{k_{0},m}\|f|L_{r}(Q_{k_{0},m})\|^{p}.
\end{equation}

Next, it is clear that
\begin{multline}
\label{eq4.20}
\Bigl(\sum\limits_{k=k_{0}}^{\infty}\Bigl(\sum\limits_{m \in \mathbb{Z}^{n}} t^{p}_{k,m} \Bigl(\mathcal{E}_{l}(\widetilde{f}, \widetilde{c}Q_{k,m})_{r}\Bigr)^{p}\Bigr)^{\frac{q}{p}}\Bigr)^{\frac{1}{q}}
\le  \\
\le
C \sum\limits_{i=1}^{3}\Bigl(\sum\limits_{k=k_{0}}^{\infty}\Bigl(\sum\limits_{m \in A_{i}(k)} t^{p}_{k,m} \Bigl(\mathcal{E}_{l}(\widetilde{f}, \widetilde{c}Q_{k,m})_{r}\Bigr)^{p}\Bigr)^{\frac{q}{p}}\Bigr)^{\frac{1}{q}}
=:S_{1}+S_{2}+S_{3}.
\end{multline}

To estimate $S_{1}$, we employ inequality \eqref{eq4.9} and Theorem~\ref{Th2.1}, while an estimate of $S_{3}$ will depend on
inequalities \eqref{eq4.16},\eqref{eq4.19}  and Theorem~\ref{Th2.1}. The sum $S_{2}$ is estimated in a~clear way.

As a~result, we have
\begin{equation}
\begin{split}
\label{eq4.21}
&\Bigl(\sum\limits_{k=k_{0}}^{\infty}\Bigl(\sum\limits_{m \in \mathbb{Z}^{n}} t^{p}_{k,m} \Bigl(\mathcal{E}_{l}(\widetilde{f}, \widetilde{c}Q_{k,m})_{r}\Bigr)^{p}\Bigr)^{\frac{q}{p}}\Bigr)^{\frac{1}{q}} \le\\
& \le C \Bigl(\sum\limits_{k=k_{0}}^{\infty}\Bigl(\sum\limits_{\substack{m \in \mathbb{Z}^{n} \\ Q^{\ast}_{k,m} \subset \Omega}}t^{p}_{k,m}\Bigl(\mathcal{E}_{l}(f,Q^{\ast}_{k,m})_{r}\Bigr)^{p}\Bigr)^{\frac{q}{p}}\Bigr)^{\frac{1}{q}}+C \sum\limits_{\substack{m \in \mathbb{Z}^{n} \\ Q_{k_{0},m} \subset \Omega}}t^{p}_{k_{0},m}\|f|L_{r}(Q_{k_{0},m})\|^{p}.
\end{split}
\end{equation}

Next, from \eqref{eq2.6}, \eqref{eq2.13} and \eqref{eq4.21} it is easy to see that
\begin{equation}
\begin{split}
\label{eq4.22}
&\Bigl(\sum\limits_{k=1}^{\infty}\Bigl(\sum\limits_{m \in \mathbb{Z}^{n}} t^{p}_{k,m} \Bigl(\mathcal{E}_{l}(\widetilde{f}, \widetilde{c}Q_{k,m})_{r}\Bigr)^{p}\Bigr)^{\frac{q}{p}}\Bigr)^{\frac{1}{q}} \le
 C \Bigl(\sum\limits_{k=k_{0}}^{\infty}\Bigl(\sum\limits_{m \in \mathbb{Z}^{n}} t^{p}_{k,m} \Bigl(\mathcal{E}_{l}(\widetilde{f}, \widetilde{c}Q_{k,m})_{r}\Bigr)^{p}\Bigr)^{\frac{q}{p}}\Bigr)^{\frac{1}{q}} + \\
&+ C \sum\limits_{m \in \mathbb{Z}^{n}}t^{p}_{0,m}\|\widetilde{f}|L_{r}(Q_{0,m})\|^{p}
\end{split}
\end{equation}

Note that the constant $C > 0$ in \eqref{eq4.22} depends on $k_{0}$.

Now from \eqref{eq4.21}, \eqref{eq4.22} it follows that the proof of the theorem will be completed once we estimate
$$
\Sigma:=\sum\limits_{\substack{m \in \mathbb{Z}^{n} \\ Q_{0,m} \bigcap \mathbb{R}^{n} \setminus \overline{\Omega} \neq \emptyset}}t^{p}_{0,m}\|\widetilde{f}|L_{r}(Q_{0,m})\|^{p}.
$$

We fix an index $m \in \mathbb{Z}^{n}$ such that $Q_{0,m} \bigcap \mathbb{R}^{n} \setminus \overline{\Omega} \neq \emptyset$. Using \eqref{eq4.3}, this gives (we also use simple estimate $\|P_{Q}[f]|L_{r}(Q)\| \le C E_{l}(f,Q)_{r} + C \|f|L_{r}(Q)\|$ which holds true for near best approximant $P_{Q}[f]$)
\begin{multline}\label{eq4.23}
t^{p}_{0,m}\Bigl(\|\widetilde{f}|L_{r}(Q_{0,m})\|^{r}\Bigr)^{\frac{p}{r}} \le C t^{p}_{0,m}\Bigl(\sum\limits_{\substack{Q \in F \\ Q \subset c Q_{0,m} \bigcap \Omega}}\Bigl(E_{l}(f,Q)_{r}\Bigr)^{r}+\|f|L_{r}(Q)\|^{r} + \|f|L_{r}(Q_{0,m} \bigcap \Omega)\|\Bigr)^{\frac{p}{r}} \le \\
 \le C t^{p}_{0,m}\|f|L_{r}(c Q_{0,m} \bigcap \Omega)\|^{p} + C S_{m}.
\end{multline}

The constant $c \geq 1$ in \eqref{eq4.23} depends only on $l, \varepsilon, \delta$ and $n$.

Arguing as in the derivation of estimate \eqref{eq4.9}, we have
\begin{equation}
\label{eq4.24}
\Bigl(\sum\limits_{\substack{m \in \mathbb{Z}^{n} \\ Q_{0,m} \bigcap \mathbb{R}^{n} \setminus \overline{\Omega} \neq \emptyset}}S_{m}\Bigr)^{\frac{1}{p}} \le C \Bigl(\sum\limits_{j=1}^{\infty}2^{-j\mu\alpha_{1}}\Bigl(\sum\limits_{\substack{ m \in  \mathbb{Z}^{n} \\ Q_{j,m} \subset C Q_{0,m} \bigcap \Omega}}t^{p}_{j,m}\Bigl(\mathcal{E}_{l}(f,Q^{\ast}_{j,m})_{r}\Bigr)^{p}\Bigr)^{\frac{\mu}{p}}\Bigr)^{\frac{q}{\mu}\frac{1}{q}}.
\end{equation}

Applying the H\"older inequality  with exponents $\frac{q}{\mu}$, $\Bigl(\frac{q}{\mu}\Bigr)'$, to the right-hand side of \eqref{eq4.24} and using finite overlap multiplicity (which depends only on $\varepsilon$, $\delta$, $n$, $l$) of cubes $cQ_{0,m}$, this gives in
combination with \eqref{eq4.23}
\begin{equation}
\label{eq4.25}
\begin{split}
&\Bigl(\sum\limits_{m \in \mathbb{Z}^{n}}t^{p}_{0,m}\Bigl(\|\widetilde{f}|L_{r}(Q_{0,m})\|^{r}\Bigr)^{\frac{p}{r}}\Bigr)^{\frac{1}{p}} \le \\
&\le C \Bigl(\sum\limits_{k=1}^{\infty}\Bigl(\sum\limits_{\substack{m \in \mathbb{Z}^{n} \\ Q^{\ast}_{k,m} \subset \Omega}}t^{p}_{k,m}\Bigl(\mathcal{E}_{l}(f,Q^{\ast}_{k,m})_{r}\Bigr)^{p}\Bigr)^{\frac{q}{p}}\Bigr)^{\frac{1}{q}}+C \Bigl(\sum\limits_{\substack{m \in \mathbb{Z}^{n} \\ Q_{0,m} \subset \Omega }}t^{p}_{0,m}\left\|f|L_{r}\Bigl(Q_{0,m}\bigcap \Omega\Bigr)\right\|^{p}\Bigr)^{\frac{1}{p}}.
\end{split}
\end{equation}

Now \eqref{eq4.6} follows from \eqref{eq4.21}, \eqref{eq4.22} and \eqref{eq4.25}. This completes the proof of the theorem.

\end{document}